\documentclass[reqno]{amsart}

\usepackage[all]{xy}
\usepackage{graphicx}

\usepackage{amssymb}
\usepackage{amsmath}
\usepackage{mathrsfs}
\usepackage{epsfig}
\usepackage{amscd}
\usepackage{graphicx, color}

\def\E{\ifmmode{\mathbb E}\else{$\mathbb E$}\fi} 
\def\N{\ifmmode{\mathbb N}\else{$\mathbb N$}\fi} 
\def\R{\ifmmode{\mathbb R}\else{$\mathbb R$}\fi} 
\def\Q{\ifmmode{\mathbb Q}\else{$\mathbb Q$}\fi} 
\def\C{\ifmmode{\mathbb C}\else{$\mathbb C$}\fi} 
\def\H{\ifmmode{\mathbb H}\else{$\mathbb H$}\fi} 
\def\Z{\ifmmode{\mathbb Z}\else{$\mathbb Z$}\fi} 
\def\P{\ifmmode{\mathbb P}\else{$\mathbb P$}\fi} 
\def\T{\ifmmode{\mathbb T}\else{$\mathbb T$}\fi} 
\def\SS{\ifmmode{\mathbb S}\else{$\mathbb S$}\fi} 
\def\DD{\ifmmode{\mathbb D}\else{$\mathbb D$}\fi} 

\newcommand{\del}{\partial}

\newcommand{\ben}{\begin{enumerate}}
\newcommand{\een}{\end{enumerate}}
\newcommand{\be}{\begin{equation}}
\newcommand{\ee}{\end{equation}}
\newcommand{\bea}{\begin{eqnarray}}
\newcommand{\eea}{\end{eqnarray}}
\newcommand{\beastar}{\begin{eqnarray*}}
\newcommand{\eeastar}{\end{eqnarray*}}
\newcommand{\bc}{\begin{center}}
\newcommand{\ec}{\end{center}}

\theoremstyle{theorem}
\newtheorem{thm}{Theorem}[section]
\newtheorem{cor}[thm]{Corollary}
\newtheorem{lem}[thm]{Lemma}
\newtheorem{prop}[thm]{Proposition}

\theoremstyle{definition}
\newtheorem{defn}[thm]{Definition}
\newtheorem{rem}[thm]{Remark}

\newtheorem{notation}[thm]{\rm\bfseries{Notation}}

\newtheorem*{thm*}{Theorem}

\numberwithin{equation}{section}

\def\R{{\mathbb R}}

\def\E{{\mathbb E}}
\def\Z{{\mathbb Z}}
\def\C{{\mathbb C}}
\def\R{{\mathbb R}}
\def\P{{\mathbb P}}

\def\N{{\mathbb N}}

\def\11{{\mathbb I}}

\def\delbar{{\overline \partial}}

\def\C{\mathbb{C}}
\def\Z{\mathbb{Z}}

\def\T{\mathbb{T}}

\def\Q{\mathbb{Q}}

\def\E{\ifmmode{\mathbb E}\else{$\mathbb E$}\fi} 
\def\N{\ifmmode{\mathbb N}\else{$\mathbb N$}\fi} 
\def\R{\ifmmode{\mathbb R}\else{$\mathbb R$}\fi} 
\def\Q{\ifmmode{\mathbb Q}\else{$\mathbb Q$}\fi} 
\def\C{\ifmmode{\mathbb C}\else{$\mathbb C$}\fi} 
\def\H{\ifmmode{\mathbb H}\else{$\mathbb H$}\fi} 
\def\Z{\ifmmode{\mathbb Z}\else{$\mathbb Z$}\fi} 
\def\P{\ifmmode{\mathbb P}\else{$\mathbb P$}\fi} 
\def\SS{\ifmmode{\mathbb S}\else{$\mathbb S$}\fi} 
\def\DD{\ifmmode{\mathbb D}\else{$\mathbb D$}\fi} 

\def\R{{\mathbb R}}

\def\E{{\mathbb E}}
\def\Z{{\mathbb Z}}
\def\C{{\mathbb C}}
\def\R{{\mathbb R}}

\def\N{{\mathbb N}}

\def\delbar{{\overline \partial}}
\def\id{{\operatorname{id}}}

\def\CD{{\mathcal D}}

\def\CL{{\mathcal L}}

\def\darr#1{\raise1.5ex\hbox{$\leftrightarrow$}
\mkern-16.5mu #1}

\def\roughly#1{\raise.3ex\hbox{$#1$\kern-.75em
\lower1ex\hbox{$\sim$}}}

\def\opname#1{\mathop{\kern0pt{\rm #1}}\nolimits}
\def\Re{\opname{Re}}
\def\Im{\opname{Im}}

\def\lcs{\mathfrak{lcs}}

\DeclareFontFamily{U}{MnSymbolC}{}
\DeclareSymbolFont{MnSyC}{U}{MnSymbolC}{m}{n}
\DeclareFontShape{U}{MnSymbolC}{m}{n}{
    <-6>  MnSymbolC5
   <6-7>  MnSymbolC6
   <7-8>  MnSymbolC7
   <8-9>  MnSymbolC8
   <9-10> MnSymbolC9
  <10-12> MnSymbolC10
  <12->   MnSymbolC12}{}
\DeclareMathSymbol{\intprod}{\mathbin}{MnSyC}{'270}

\begin{document}
\quad \vskip1.375truein

\def\mq{\mathfrak{q}}
\def\mp{\mathfrak{p}}
\def\mH{\mathfrak{H}}
\def\mh{\mathfrak{h}}
\def\ma{\mathfrak{a}}
\def\ms{\mathfrak{s}}
\def\mm{\mathfrak{m}}
\def\mn{\mathfrak{n}}
\def\mz{\mathfrak{z}}
\def\mw{\mathfrak{w}}
\def\Hoch{{\tt Hoch}}
\def\mt{\mathfrak{t}}
\def\ml{\mathfrak{l}}
\def\mT{\mathfrak{T}}
\def\mL{\mathfrak{L}}
\def\mg{\mathfrak{g}}
\def\md{\mathfrak{d}}
\def\mr{\mathfrak{r}}

\title[Contact triad connection]{Contact triad connection of contact manifolds
in almost contact moving frame}

\author{Yong-Geun Oh}
\address{Center for Geometry and Physics, Institute for Basic Science (IBS),
77 Cheongam-ro, Nam-gu, Pohang-si, Gyeongsangbuk-do, Korea 790-784
\& POSTECH, Gyeongsangbuk-do, Korea}
\email{yongoh1@postech.ac.kr}
\thanks{The author is supported by the IBS project \# IBS-R003-D1}
%

\begin{abstract} The notion of \emph{contact triad connection} on
contact triads $(M,\lambda,J)$ was introduced by Wang and
the present author  in early 2010's \emph{from scratch} as the contact analog to the canonical connection 
of an almost K\"ahler manifolds. The connection facilitates the study of analysis of
the contact instanton equation which ranges from local elliptic priori estimates,
for both the interior \cite{oh-wang:CR-map1,oh-wang:CR-map2}
and the boundary \cite{oh-yso:index}, to the generic perturbation theory of the asymptotic operators \cite{kim-oh:asymp-analysis},
which also encompasses the case of pseudoholomorphic curves on noncompact symplectic manifolds 
with cylindrical ends \cite{oh-kim:matrix}. The main purpose of the present paper is to
give a simpler and more canonical construction of the contact triad connection by first
giving its characterization in terms of the associated \emph{almost contact structure} and then providing its construction 
using the \emph{almost contact moving frame} in the framework of (metric) contact geometry.
\end{abstract}



\maketitle

\tableofcontents

\section{Introduction}
\label{sec:intro}

A contact manifold $ (M,\xi)$ is a $2n+1$ dimensional manifold
equipped with a completely non-integrable distribution of rank $2n$,
called a contact structure. Complete non-integrability of $\xi$
can be expressed by the non-vanishing property
$$
\lambda \wedge (d\lambda)^n \neq 0
$$
for a one-form $\lambda$ which
defines the distribution, i.e., $\ker \lambda = \xi$. Such a
one-form $\lambda$ is called a contact form associated to $\xi$.
Associated to the given contact form $\lambda$, we have the unique vector field $R_\lambda$
determined by
$$
R_\lambda \intprod \lambda \equiv 1, \quad R_\lambda \intprod d\lambda \equiv 0.
$$
The choice of $R_\lambda$ also induces the splitting 
$TM = \xi \oplus \R \langle R_\lambda \rangle$ of the canonical short exact sequence
\be\label{eq:shortexact1}
0 \to \xi \to TM \to TM/\xi \to 0,
\ee
and the splitting
$$
T^*M = (TM/\xi)^\perp \oplus (\xi)^\perp \cong \xi^* \oplus \R \{\lambda\}
$$
of the dual sequence
\be\label{eq:shortexact2}
0 \to (TM/\xi)^\perp \to T^*M \to (\xi)^\perp \to 0
\ee
where $(\cdot)^\perp$ is the annihilator of $(\cdot)$.

In relation to the study of pseudoholomorphic curves, one consider an
endomorphism $J: \xi \to \xi$ with $J^2 = -id$ and regard $(\xi,J)$
a complex vector bundle. In the presence of the contact form $\lambda$,
one usually consider the set of $J$ that is compatible to $d\lambda$ in
the sense that the bilinear form $g_\xi = d\lambda(\cdot, J\cdot)$ defines
a Hermitian vector bundle $(\xi,J, g_\xi)$ on $Q$. We call the triple $ (M,\lambda,J)$
a \emph{contact triad}.

\subsection{Advent of contact triad connection and its axioms}

For each given contact triad, we equip $Q$ with the triad metric
$$
g = d\lambda(\cdot, J \cdot) + \lambda \otimes \lambda.
$$
For a given Riemann surface $(\Sigma,j)$ and a smooth map $w:\Sigma \to M$, we
decompose the derivative $dw: T\Sigma \to TQ$ into
$$
dw  =  \pi\, dw + w^*\lambda\, R_\lambda
$$
as a $w^*TQ$-valued one-form and then
$$
d^\pi w: = \pi dw = \del^\pi w + \delbar^\pi w
$$
as a $w^*\xi$-valued one-form on $\Sigma$ where $\del^\pi w$ (resp. $\delbar^\pi w$) is the
$J$-linear part (resp. anti-$J$-linear part) of $d^\pi w$.

\begin{defn}[Contact Cauchy-Riemann map]
Let $(\Sigma, j)$ be a Riemann surface with a finite number of marked points and denote by $\dot \Sigma$
the associated punctured Riemann surfaces. We call a map $w:\dot \Sigma \to M$ a
\emph{contact Cauchy-Riemann map} if $\delbar^\pi w = 0$.
\end{defn}

It turns out that to establish the geometric analysis necessary for the
study of associated moduli space, one needs to augment the equation $\delbar^\pi w = 0$ by
\be\label{eq:closed}
d(w^*\lambda \circ j) = 0.
\ee
We refer to \cite{hofer:gafa} for the origin of such equation in
contact geometry and to \cite{oh-wang:CR-map1}, \cite{oh-wang:CR-map2}  for the detailed analytic study of priori
$W^{k,2}$-estimates and asymptotic convergence on punctured Riemann surfaces.

\begin{defn}[Contact instanton] Let $\Sigma$ be as above.
We call a pair of $(j,w)$ of a complex structure $j$ on $\Sigma$ and a map $w:\dot \Sigma \to M$ a
\emph{contact instanton} if it satisfies
\be\label{eq:contact-instanton}
\delbar^\pi w = 0, \, \quad d(w^*\lambda \circ j) = 0.
\ee
\end{defn}

In the course of studying geometric analysis of such a map \cite{oh-wang:CR-map1,oh-wang:CR-map2,oh-yso:index},
we needed to simplify the tensorial calculations by choosing a special
connection as in (almost) complex geometry. For this purpose, a special connection, called
the \emph{contact triad connection}, of the triad $ (M,\lambda,J)$ was exploited which
had been introduced by Wang the present author in \cite{oh-wang:connection}: This is
the contact analog to Ehresmann-Libermann's notion of \emph{canonical connection}
on the almost K\"ahler manifold $(M,\omega,J)$.
(See \cite{ehres-liber,libermann54,libermann55}, \cite{kobayashi:canonical}, \cite{gaudu} for general exposition
on the canonical connection.)

The following is the list of axioms that characterizes the contact triad connection.
\begin{thm}[Contact triad connection, \cite{oh-wang:connection}]
\label{thm:mainI} Let $ (M,\lambda,J)$ be any contact triad of
contact manifold $ (M,\xi)$. Denote by
$$
g_\xi + \lambda \otimes \lambda =: g
$$
the natural Riemannian metric on $Q$ induced by $(\lambda,J)$, which we call a
contact triad metric. Then there exists a unique affine connection $\nabla$ that has the
following properties:
\begin{enumerate}
\item $\nabla$ is a Riemannian connection of the triad metric.
\item The torsion tensor of $\nabla$ satisfies $T(R_\lambda, Y)=0$ for all $Y \in TQ$.
\item $\nabla_{R_\lambda} R_\lambda = 0$ and $\nabla_Y R_\lambda\in \xi$, for $Y\in \xi$.
\item $\nabla^\pi := \pi \nabla|_\xi$ defines a Hermitian connection of the vector bundle
$\xi \to M$ with Hermitian structure $(d\lambda, J)$.
\item The $\xi$ projection, denoted by $T^\pi: = \pi T$, of the torsion $T$
is of $(0,2)$-type in its complexification, i.e., satisfies the following properties:
\beastar\label{eq:TJYYxi}
T^\pi(JY,Y) = 0
\eeastar
for all $Y$ tangent to $\xi$.
\item For $Y\in \xi$,
$$
\nabla_{JY} R_\lambda+J\nabla_Y R_\lambda= 0.
$$
\end{enumerate}
We call $\nabla$ the contact triad connection.
\end{thm}

We recall that the leaf space of Reeb foliations of the contact triad $ (M,\lambda,J)$
canonically carries a (non-Hausdorff) almost K\"ahler structure which we denote by
$
(\widehat M,\widehat{d\lambda},\widehat J).
$
We would like to note that Conditions (4) and (5)
are nothing but properties of the canonical connection
on the tangent bundle of the (non-Hausdorff) almost
K\"ahler manifold $(\widehat M, \widehat{d\lambda},\widehat J_\xi)$
lifted to $\xi$. On the other hand, Conditions (1), (2), (3)
indicate this connection behaves like the Levi-Civita connection when the Reeb direction $R_\lambda$ get involved.
\emph{Condition (6) is an extra condition to connect the information in $\xi$ part and $R_\lambda$ part},
which is used to dramatically simplify our calculation in \cite{oh-wang:CR-map1}, \cite{oh-wang:CR-map2}.
The condition is the most novel part of the axiom, \emph{which has no prescient in other literature}.

Coining Axiom (6) in advance would have been unimaginable for us without the guidance of our efforts
to simplify the complicated tensorial calculations, \emph{when we had used the more common Levi-Civita
connection}, which arises in the analytic study of  
contact Cauchy-Riemann maps carried  out in \cite{oh-wang:CR-map1}. 
In hindsight,  the axiom turns out to become 
more canonical to introduce \emph{in terms of an almost contact frame} in that
it carries a natural
representation written as in Proposition \ref{prop:axiom-6} below in (almost) Cauchy-Riemann geometry.

\begin{rem}
The contact triad connection $\nabla$ canonically induces
a Hermitian connection for the Hermitian vector bundle $(\xi,J,g_\xi)$
with $g_\xi = d\lambda(\cdot, J \cdot)|_\xi$. 
We call this vector bundle connection the \emph{contact Hermitian connection}.
In fact, the triad connection 
can be generalized to the setting of \emph{almost contact manifolds} \cite{gray}.
See Remark \ref{rem:ac-extension} for further relevant remarks.
\end{rem}

In fact, one can explicitly write the formula for the contact triad connection in
terms of the Levi-Civita connection

\begin{thm}[\cite{oh-wang:connection}]\label{thm:triad-LC} Let $ (M,\lambda, J)$ be a contact triad and
$\Pi: TQ \to TQ$ be the idempotent with $\Im \Pi = \xi$ and $\ker \Pi = \R\{R_\lambda\}$.
Extend $J$ to $TQ$ setting $J(R_\lambda) = 0$ and still denote by the same $J$.
Let $\nabla$ be the contact triad connection and $\nabla^{LC}$
be the Levi-Civita connection of the triad metric. Then we have
$$
\nabla = \nabla^{LC} + B
$$
where $B$ is the tensor of type $(\begin{matrix}1\\2 \end{matrix})$ defined by
$$
B = B_1 + B_2 + B_3
$$
where
\bea\label{eq:B}
B_1(Z_1,Z_2) & = & - \frac{1}{4}\big((\nabla^{LC}_{JZ_2}J) \Pi Z_1 + J(\nabla^{LC}_{\Pi Z_2}J) \Pi Z_1
+ 2J(\nabla^{LC}_{\Pi Z_1}J)\Pi Z_2\big),\nonumber\\
B_2(Z_1,Z_2) & = & \frac{1}{8}\big(-\langle J(\CL_{R_\lambda} J)Z_1, Z_2\rangle
+\langle JZ_1, Z_2\rangle\big) R_\lambda,
\nonumber\\
B_3(Z_1,Z_2) & = &\frac{1}{2}\big( -\langle Z_2, R_\lambda\rangle JZ_1
-\langle Z_1, R_\lambda\rangle JZ_2
+ \langle JZ_1, Z_2\rangle R_\lambda\big).
\eea
\end{thm}
Complication of this explicit formula of the contact triad connection relative to the Levi-Civita connection
shows that it is highly nontrivial to guess the existence of such a connection, let alone its
uniqueness. 

\subsection{Almost contact characterization of contact triad connection}

While our introduction of Condition (6) is motivated by our attempt to simplify the tensor
calculations \cite{oh-wang:CR-map1}, it has a nice geometric interpretation in terms of
$CR$-geometry.  We call a one-form $\alpha$ is $CR$-holomorphic
if $\alpha$ satisfies
\beastar
\nabla_{R_\lambda} \alpha & = & 0, \\
\nabla_Y \alpha + J \nabla_{JY} \alpha & = & 0 \quad \text{for }\, Y \in \xi.
\eeastar

\begin{prop} In the presence of other defining
properties of contact triad connection, condition (6) is equivalent to the statement that
$\lambda$ is $J$ holomorphic in the $CR$-sense.
\end{prop}

An immediate consequence of this proposition is the following
important property of the pull-back form $w^*\lambda$ for any contact Cauchy-Riemann map
which plays an important role in the study of asymptotic behaviors of contact instantons
\cite{oh:contacton,oh-yso:index}, i.e., of solutions
$w: \dot \Sigma \to M$  of \eqref{eq:contact-instanton} near punctures of $\dot \Sigma$.

\begin{cor} Let $(\Sigma,j)$ be a compact Riemann surface and $\dot \Sigma$ be a
punctured surface arising by the removal of a finite number of points from $\Sigma$.
Let $w: (\dot \Sigma,j) \to  (M,J)$ be any Cauchy-Riemann map. Then
 the complex valued one-form 
 $$
 \chi: = w^*\lambda \circ j + \sqrt{-1} w^*\lambda
 $$
has a decomposition $\chi = h_{(w,j)} + df_{(w,j)}$ where 
$h_{(w,j)}$ is an abelian differential of the Riemann surface $(\dot \Sigma, j)$.
\end{cor}
Some further elaboration of this point of view in relation to the study of compactified moduli
space of higher genus contact instantons is given in \cite{oh:contacton}, \cite{oh:abelian-differential}.

This motivates us to re-visit  Wang and the present author's earlier construction 
\cite{oh-wang:connection} of contact triad connection, which was originally made from scratch, 
by adapting Kobayashi's complex moving frame method  \cite{kobayashi:canonical} 
applied for the study of structure equations of almost Hermitian
connection in the framework of \emph{almost contact manifolds}. The outcome is a new construction of
contact triad connections along a simpler and more canonical route an outline of which is
now in order.

We first reformulate the defining conditions of the triad connection in terms of
the associated almost contact manifold and the reduction of the structure group $SO(2n+1)$ to
$U(n) \times \{1\}$. (See \cite{gray,casals-pancho-presas,borman-eliash-murphy}
for the definition and its usage in relation to
the existence question of contact structure.)
Choose a local orthogonal Darboux  frame of $TM = \xi \oplus\R \langle R_\lambda \rangle $ given by
$$
 E_1,\cdots, E_n, F_1,\cdots, F_n,  R_\lambda, \quad F_i: = J E_i
 $$
and denote its dual co-frame by
$
e^1,\cdots, e^n, f^1,\cdots, f^n, \, \lambda.
$ 
where $\xi^{(1,0)}$ is the $(1,0)$ part of $\xi_\C = (\xi,J) \otimes \C$. 
Consider the complexification of $TM$ and complex linear extension $g_\C = g \otimes \C$
of the contact triad metric  $g$.

\begin{defn}[Almost contact moving frame] Let $(M,\lambda, J)$ be a contact triad equipped with triad metric.
 We call a \emph{almost contact frame}
the associated complexification of a Darboux frame
$$
\C\{\eta_1,\cdots, \eta_n, \overline \eta_1, \cdots, \overline \eta_n\} \oplus \C \{R_\lambda\}, \quad \eta_i : = \frac1{\sqrt{2}}(E_i- \sqrt{-1}  F_i).
$$
\end{defn}
We consider its dual almost contact frame
$$
\{\theta^1, \cdots, \theta^n, \overline \theta^1, \cdots, \overline \theta^ n\} \oplus \C \langle \lambda \rangle
$$
on the complex vector bundle with decomposition
\be\label{eq:ac-decomposition}
 (\xi^* \otimes \C) \oplus \C \langle \lambda \rangle = 
 \xi^{*(1,0)} \oplus \xi^{*(0,1)} \oplus \C \langle \lambda \rangle.
\ee
We now express a connection on the complex subbundle 
$$
\xi^{(1,0)} \oplus \C \langle R_\lambda \rangle \subset TM \otimes \C
$$
in terms of the almost contact frame.

 Then $g_\C$ restricts to a Hermitian
metric on the complex vector bundle 
$\xi^{(1,0)} \oplus \{0\} \subset \xi^{(1,0)} \oplus \C\langle R_\lambda \rangle$.
Its dual frame is given by
$$
\{\theta^1, \cdots, \theta^n\} \oplus \{\lambda\}
$$
where we recall that $\theta^i$ is 1-forms of $(1,0)$-type and $\lambda$ is a real one-form.
The connection form of the complex vector bundle $\xi^{(1,0)} \oplus \C \langle R_\lambda \rangle \to M$
associated the above frame is a matrix-valued one-form  consisting of a skew-Hermitian matrix
$$
\omega_j^i, \quad i, \, j = 1,\cdots, n
$$
and of complex (column) vectors
$$
\alpha^0 = (\alpha^0_j) \quad \beta_0^t = (\beta_0^i) = (\beta^i), \quad i,\, j = 0, \cdots, n
$$
of one-forms. We will first show that the vanishing
$$
\alpha_0^0= 0
$$
is equivalent to Axiom (3).

We also have the following translation of Axiom (6) above, which is a simplified version of
Proposition \ref{prop:axiom-3} in the main text. 
We refer to Proposition \ref{prop:axiom-3} for further generalization.
In hindsight, this axiom, combined with other axioms, seems to
explain the `dramatic simplification' of the tensor calculations occurring in the study of
geometric analysis of contact instantons and  pseudoholomorphic curves on the symplectization,
(See \cite{oh-wang:CR-map1,oh-wang:CR-map2}, \cite{oh-yso:index}, 
\cite{kim-oh:asymp-analysis} for the analysis of contact instantons, and \cite{oh-kim:matrix} for 
the analysis of pseudoholomorphic curves via that of contact instantons.)

\begin{prop}[Compare with Proposition \ref{prop:axiom-3}]\label{prop:axiom-6} 
Assume $c = 0$. We decompose 
$$
\alpha_k^0 = (\alpha_k^{0\pi}) + \alpha_k^{0\perp} =  (\alpha_k^{0\pi})^{(1,0)} + 
 (\alpha_k^{0\pi})^{(0,1)} + \alpha_k^{0\perp}
 $$
 according to \eqref{eq:ac-decomposition}. Then we have $\alpha_k^{0\perp}=0$ and
  Axiom (6) above is equivalent to
\be\label{eq:alpha0-antiholomorphic}
(\alpha_k^{0\pi})^{(1,0)} = 0
\ee
i.e., $\alpha^{0\pi}_k$ is anti-holomorphic.
\end{prop}

Now the above discussion enables us to translate the condition on the real torsion $T$
laid out in Definition \ref{defn:triad-connection} to those of complex torsion $\Theta$ as follows.
\begin{prop}[Compare with Proposition \ref{prop:in-movingframe}]\label{prop:in-movingframe-intro}
Let $\nabla$ be any connection $\nabla g = 0$. Equip it with an almost contact frame
and let $\Theta$ be as above. Then the following equivalence holds:
\begin{enumerate}
\item Axiom (6) is equivalent to \eqref{eq:alpha0-antiholomorphic}.
\item Axioms (2), (3) together is  equivalent to $R_\lambda \intprod \Theta^0 = 0$.
\item Conditions (1) and (4)  are equivalent to $\Theta^\pi$ is Hermitian.
\item Condition (5) is equivalent to the vanishing $(\Theta^\pi)^{(1,1)} = 0$.
\end{enumerate}
\end{prop}

This completes translation of the defining conditions given in Definition \ref{defn:triad-connection}
into the language of the almost contact frame.

Then the following simple complex characterization of the contact triad connection from Theorem \ref{thm:mainI}
and the proof of its unique existence is the main results of the present paper. 

\begin{thm}[Main Theorem] \label{thm:main} Assume $c = 0$.
Let $\Theta$ be the torsion form 
with respect to an almost contact frame.
We denote by $\Theta^\pi: = \Pi^*\Theta$ as a $\xi$-valued 2-form. There exists a unique connection 
associated to each contact triad $(M,\lambda,J)$ such that
its (complexified) torsion form $\Theta$ with respect to the frame satisfies 
\be\label{eq:Theta0}
\Theta^0 = d\lambda + \sum_{k=1}^n \alpha_k^0 \wedge \theta^k, \quad \Theta^{\pi(1,1)} = 0
\ee
such that $\alpha_k^{0\pi}$ is  of  $(0,1)$-type, i.e., anti-holomorphic.

In addition, we also have $\Theta^{\pi(2,0)} = 0$ and so $\Theta^\pi$ of $(0,2)$-type.
\end{thm}
The whole section of 
Section \ref{sec:construction} will be occupied by the proof of the generalized version thereof
considering the general case of \emph{$c$-contact triad connection} with $c \neq 0$. We refer thereto
for the more general statement in this general context.

We can also canonically lift the connection to
the $\lcs$-fication $(M \times S^1, \omega^{\mathfrak b})$ of contact manifold $(M,\lambda)$ with
$$
\omega^{\mathfrak b} = d\lambda + \mathfrak{b}_R \wedge \lambda, \quad \mathfrak{b}_R: = \frac1R d\theta
$$
which we call the $\lcs$-fications of a triad connection. This also covers the case
$$
(M \times \R, d \lambda + ds \wedge \lambda)
$$
as the `zero-temperature limit'  thereof. This enters in the analysis of pseudoholomorphic curves on
the $\lcs$-fication including the case of the symplectization of contact manifold $(M,\lambda)$.
We refer readers to \cite{oh-savelyev}, \cite{oh:book-contact} for a more elaborate discussion on this lifting.

Organization of the paper is now in order. In Section \ref{sec:can-Kahler}, we recall the basic
properties of the canonical connection of almost K\"ahler manifolds collected from \cite{kobayashi-nomizu2},
\cite{kobayashi:canonical}, \cite{oh:book1} to motivate the construction of contact triad connection.
Section \ref{sec:triad-connection} reviews the axioms of contact triad connection from \cite{oh-wang:connection}
and collect its general properties which will be crucial for the complex translation of 
contact triad connection. These have been scattered around in the author's earlier articles
\cite{oh-wang:CR-map1,oh-wang:CR-map2}, \cite{oh:contacton}, \cite{oh-yso:index}.
Section \ref{sec:in-complexframe} provides the complex translation of the axioms of 
the contact triad connection and Proposition \ref{prop:in-movingframe-intro}. Then 
we prove the main result Theorem \ref{thm:main} in Section \ref{sec:construction}. 
\bigskip

\noindent{\bf Acknowledgement:} We thank Rui Wang for our collaboration of 
the construction of contact triad connection in \cite{oh-wang:connection},
and  useful comments on the early version of the present paper. We sincerely thank the unknown
referee for her/his careful reading of the original manuscript and pointing out many careless 
errors in tensor calculations. Accommodation of the referees' comments and other elaboration of
some omitted details in this revision have improved the exposition which we hope makes the paper
more readable.

\section{Canonical connection of almost K\"ahler manifolds}
\label{sec:can-Kahler}

We start with outlining Liebermann's construction \cite{libermann55} of canonical connection on almost
Hermitian manifolds. To connect his construction to Chern connection  \cite{chern:connection}
in the integrable case, it is best to describe it in terms of its complexification following the exposition of
Kobayashi \cite{kobayashi:canonical}.

Let $(M,J)$ be an almost complex manifold, i.e., $J$ be a bundle
isomorphism $J : TM \to TM$ over the identity satisfying $J^2 = -
\id$.
\begin{defn} A metric $g$ on $(M,J)$ is called \emph{(almost) Hermitian}, if $g$ satisfies
$$
g(Ju,Jv) = g(u,v), \quad u, \, v \in T_x M, \, x \in M.
$$
We call the triple $(M,J,g)$ an almost Hermitian manifold.
\end{defn}

For any given almost Hermitian manifold $(M,J,g)$, the bilinear form
\be\label{eq:Phi}
\Phi : = g(J \cdot, \cdot)
\ee
is called the \emph{fundamental two-form} in \cite{kobayashi-nomizu2}, which
is nondegenerate.

\begin{defn} An almost Hermitian manifold $(M,J,g)$ is an \emph{almost
K\"ahler manifold} if the two-form $\Phi$ above is closed.
\end{defn}

\begin{defn}
A (almost) Hermitian connection $\nabla$ is an affine connection
satisfying
$$
\nabla g = 0 = \nabla J.
$$
\end{defn}
Such a connection always exists.

In general the torsion $T = T_\nabla$ of the almost Hermitian connection $\nabla$ is not zero, 
even when $J$ is integrable. By definition, $T$ is a $(2,1)$-tensor defined by
$$
T(X,Y) = \nabla_X Y - \nabla_Y X - [X,Y].
$$
The following is the almost complex version of the Chern connection
\cite{chern:connection} in complex
geometry.

\begin{thm}[\cite{gaudu}, \cite{kobayashi:canonical}]\label{thm:canonicalconnection}
On any almost Hermitian manifold $(M,J,g)$,
there exists a unique Hermitian connection $\nabla$ on $TM$ satisfying
\be\label{eq:canonical-nabla}
T(X,JX) = 0
\ee
for all $X \in TM$.
\end{thm}

\begin{defn}\label{defn:canonical-nabla} A canonical connection of an almost Hermitian connection is
defined to be one that has the torsion property \eqref{eq:canonical-nabla}.
\end{defn}

Next we would like to derive how the \emph{almost K\"ahler} property
$d\Phi = 0$ for the bilinear form $\Phi = g(\cdot, J\cdot)$ affects
on the torsion property of the canonical connection. Such a circumstance naturally arises 
through the choice of $\omega$-compatible almost complex structure $J$ on a symplectic
manifold $(M,\omega)$, i.e., $g = \omega(\cdot, J\cdot)$ and $\Phi = \omega$.

For this purpose, it is best to
express the connection in the  complexified tangent bundle 
$TM \otimes \C$ and to extend $J$ thereto. Then we have the eigenvalue
decomposition
$$
 TM \otimes \C = T^{(1,0)}M \oplus T^{(0,1)}M
$$
where the complex subbundle $T^{(1,0)}M$ (resp. $T^{(0,1)}M$) 
of $TM \otimes \C$ is the eigen-space of $J$ with eigenvalue $\sqrt{-1}$
(resp. $-\sqrt{-1}$).  We denote by $\Pi'$ (resp. $\Pi''$) the projection
$TM \otimes \C  \to T^{(1,0)}M$ (resp. $TM \otimes \C  \to T^{(0,1)}M$).

In terms of this complexification, the defining property 
\eqref{eq:canonical-nabla} can be expressed as follows.
Denote by $T_\C = T \otimes \C$  the complex linear
extension of the torsion tensor $T$ of the connection $\nabla$,
and $\Pi'$ is the $(1,0)$-projection to $T^{(1,0)}M$.
We call the projection
\be\label{eq:Theta}
\Theta: = \Pi' T_\C
\ee
of $T_\C$ as a $T^{(1,0)}M$-valued two-form \emph{complex torsion} of  the almost Hermitian connection. 
The $(0,2)$-component of the two-form $\Theta$ depends only on the almost complex structure $J$
and does not depend on the choice of connection, which is nothing but the complexification
(more precisely, a constant multiple of the index-lowering) of 
the Nijenjuis tensor which is the $(2,1)$-tensor $N$ defined by
\be\label{eq:nijenhuis}
N(X,Y) = [JX,JY] - [X,Y] -J[X,JY] - J[JX,Y]
\ee
for all $X,\, Y \in \mathfrak{X}(M)$.
(See \cite[Theorem 3.4]{kobayashi-nomizu2} (with different normalization factor), \cite[Equation (1.3)]{kobayashi:canonical}.)

\begin{defn} Let $(M,J,g)$ be an almost Hermitian manifold.
A Hermitian connection $\nabla$ of the complex vector bundle $T^{(1,0)}M$
is called the {\it canonical} connection, if its complex torsion 
$\Theta$ satisfies $\Theta^{(1,1)} = 0$.
\end{defn}

Straight computation shows that \eqref{eq:canonical-nabla} is equivalent to $\Theta^{(1,1)} = 0$.
Therefore Theorem \ref{thm:canonicalconnection} is equivalent  to
the following theorem.

\begin{thm}[\cite{kobayashi:canonical}]\label{thm:ac-unique} On any Hermitian manifold $(M,J,g)$, 
there exists a unique Hermitian connection $\nabla$ on $TM$ leading to the
canonical connection on $T^{(1,0)}M$. 
\end{thm}

The following is also proved by Kobayashi.

\begin{thm}[\cite{kobayashi:canonical}]\label{thm:aK-(2,0)-vanishing} Let $(M,J,g)$ be almost
K\"ahler and let $\nabla$ be the canonical connection of $T^{(1,0)}M$. Then
$\Theta^{(2,0)} = 0$ in addition, and hence $\Theta$ is of type $(0,2)$.
\end{thm}

The following properties can be derived from this latter theorem.

\begin{prop}\label{prop:TJYY} Let $(M,J, g)$ be an almost K\"ahler manifold
and $\nabla$ be its canonical connection. Denote by $T$ be its torsion tensor. Then
the following identity holds:
\be\label{eq:TJYY}
T(JY,Y) = 0
\ee
for all $Y$. Moreover
\be\label{eq:TJYTY}
T(JY,Z) = T(Y,JZ)
\ee
and
\be\label{eq:JTJYZ}
JT(JY,Z) = T(Y,Z)
\ee
for all vector fields $Y, \, Z$ on $M$.
\end{prop}
\begin{proof} For general real vector field $Z$ on $M$ can be decomposed into
$$
Z = \frac{Z - iJ Z}{2} + \frac{Z+iJZ}{2} \in T^{(1,0)}W \oplus T^{(0,1)}W.
$$
Since $\Theta = \Pi' T_\C$ is of $(0,2)$-type, to compute
$
\Theta(JY,Y)
$
we have only to consider $T^{(0,1)}$-components of $JY$ and $Y$. Therefore
\beastar
\Theta(JY,Y) & = & \Theta\left(\frac{JY+iJJY}{2},\frac{Y+iJY}{2}\right)\\
& = & \Theta\left(-i\left(\frac{Y+iJY}{2}\right),\frac{Y+iJY}{2}\right) = 0.
\eeastar
Now we compute $\Pi''T_\C$-component of $T(JY,Z)$. But since $T(JY,Y)$ is
a real vector, this component is nothing but $\overline{\Theta(JY,Y)}$
and hence also becomes zero. This finishes the proof of \eqref{eq:TJYY}.

For \eqref{eq:TJYTY}, we apply \eqref{eq:TJYY} to the vector fields
$Y + Z$ which will then gives rise to its proof.

Finally we prove \eqref{eq:JTJYZ}. Similarly as above, we compute
\beastar
\Theta(JY,Z)  = \Theta\left(\frac{JY - iY}{2},\frac{Z + iJZ}{2}\right)
= \frac{1}{2}(\Theta(JY,Z) - i\Theta(Y,Z))
\eeastar
where we use \eqref{eq:TJYTY} for the second equality. Hence $\Theta(JY,Z) = -i \Theta(Y,Z)$.
Therefore we derive
$$
T(JY,Z) = 2\Im(\Theta(Y,Z)).
$$
Then we compute
\beastar
JT(JY,Z) = 2 \Im (J \Theta(Y,Z))) = 2 \Im (i\Theta(Y,Z))) =
2 \Re(\Theta(Y,Z)) = T(Y,Z).
\eeastar
This finishes the proof of \eqref{eq:JTJYZ}.
\end{proof}

In fact, one can easily check that the properties spelled out in this proposition
completely characterize the canonical connection, which gives rise to the following
real characterization of the canonical connection.

\begin{thm}\label{thm:can-ac} Let $(M,J,g)$ be an almost K\"ahler manifold.
Then an almost Hermitian connection $\nabla$ of $(M,J,g)$ is a canonical connection if and only if
the torsion $T$ of $\nabla$ satisfies the properties \eqref{eq:TJYY}, \eqref{eq:TJYTY}
and \eqref{eq:JTJYZ}.
\end{thm}
\begin{proof} We have only to prove that any almost Hermitian connection satisfying
\eqref{eq:TJYY}, \eqref{eq:TJYTY} and \eqref{eq:JTJYZ} satisfies the $(1,1)$-part of
$\Theta= \Pi'  T_\C$ vanishes, i.e.,
$$
\Pi'  T_\C (X + iJ X, Z + iJ Z) = 0
$$
for all $X, \, Z$ on $W$. The proof of this is straightforward computation using
the given properties and so omitted.
\end{proof}

\section{Contact triad connection on contact manifolds}
\label{sec:triad-connection}

Let $(M,\xi)$ be a contact manifold and a contact form $\lambda$ of
$\xi$ be given. On $M$, the Reeb vector field $R_\lambda$ associated to the contact
form $\lambda$ is the unique vector field satisfying
\be\label{eq:Liouville} X \intprod \lambda = 1, \quad X \intprod
d\lambda = 0.
\ee
Therefore the tangent bundle $TM$ has the splitting $TM =
 \xi \oplus \R \langle R_\lambda\rangle$. We denote by
$$
\pi_\lambda: TM \to \xi
$$
the corresponding projection. We extend a complex structure $J_\xi$ on
the bundle $\xi \to M$ to the endomorphism $J: TM \to TM$ by the trivial extension of $J_\xi$
obtained by setting 
$$
J(R_\lambda): = 0.
$$
\begin{rem} In \cite[p. 47]{blair}, Blair calls such a tuple $(J,R_\lambda, \lambda, g)$ a
\emph{contact metric structure}, and \emph{$K$-contact} if it satisfies $\CL_{R_\lambda}J = 0$
in addition. (Also see Remark \ref{rem:CLRJ=0} for further related comments on this vanishing.)
\end{rem}

\subsection{Contact triads and Blair's lemma}

\begin{defn}[Contact triad] We call the triple $(M,\lambda, J)$ a \emph{contact triad}.
\end{defn}
Each contact triad $(M,\lambda, J)$  naturally induces
a metric on $M$ compatible to the contact form $\lambda$ and $J$ which is defined by
$$
g =  d\lambda(\cdot, J \cdot)  +  \lambda \otimes \lambda
$$
We call this metric the $(\lambda,J)$-compatible metric.
In this section, we associate a particular type of connection on $M$
which we call \emph{the contact triad connection} of the triad $(M,\lambda, J)$.
\begin{rem} It is important to notice that the scaling property of this metric is different depending on
the $\xi$-direction or in the Reeb direction: Relative to the contact form $\lambda$, the former scales 
linearly while the latter scales quadratically when $\lambda$ is scaled linearly.
\end{rem}

The following is a nice interplay between the triad metric and the Lie derivative of 
$J$, which we quote from Blair's monograph \cite{blair}. It plays an important role in the analysis of 
contact instanton equation in \cite{oh-wang:CR-map1}, \cite{oh-yso:index} and \cite{kim-oh:asymp-analysis}.
For the self-containedness and for the purpose of familiarizing readers with the shape of 
relevant tensor calculations throughout the entire paper,  we provide its proof here.

\begin{prop}[Lemma 6.2 \cite{blair}]\label{prop:symmetry}
For $Y,\, Z \in \xi$, both operators $\CL_{R_\lambda}J$ and $\CL_{R_\lambda}J J$ are 
symmetric with respect to the triad metric.
\end{prop}
\begin{proof} Once we prove $\CL_{R_\lambda}J$ is symmetric, then the symmetry of
$\CL_{R_\lambda}J J$ also follows by the following simple calculation: 
\beastar
(\CL_{R_\lambda}J J)^t = J^t (\CL_{R_\lambda}J)^t. = -J \CL_{R_\lambda}J
= \CL_{R_\lambda}J J
\eeastar
where the last equality follows from the vanishing $ \CL_X J J + J \CL_X J = 0$
for all $X$ which in turn arises from that $J^2 = -\Pi$ and $\Pi$ is parallel.

Now we will prove $\CL_{R_\lambda}J$ is symmetric, i.e., 
\be\label{eq:symmetry}
\langle (\CL_{R_\lambda}J)Y,Z \rangle = \langle Y, (\CL_{R_\lambda}J) Z \rangle
\ee
for all $Y, \, Z$.  If one of them is $R_\lambda$, it follows that both sides vanish
by the properties $[R_\lambda, Y] \in \xi$ for all $Y$ tangent to $\xi$ and $JR_\lambda = 0$.

Therefore we now assume $Y, \, Z \in \xi$. Using the definition of the triad metric
$g = d\lambda(\cdot, J\cdot) + \lambda \otimes \lambda$, we evaluate
$$
\langle (\CL_{R_\lambda}J)Y,Z \rangle = d\lambda \big((\CL_{R_\lambda}J)Y, JZ)
+ \lambda( (\CL_{R_\lambda}J)Y \big) \cdot \lambda(Z).
$$
We rewrite the second term into 
$$
\lambda( (\CL_{R_\lambda}J)Y) \cdot \lambda(Z)=
\big(\lambda(\CL_{R_\lambda}(JY)) - \lambda(J \CL_{R_\lambda} Y)\big) \cdot \lambda(Z)
$$
and then rewrite the first factor thereof into
\beastar
\lambda(\CL_{R_\lambda}(JY)) - \lambda(J \CL_{R_\lambda} Y) & =  & \lambda(\CL_{R_\lambda}(JY)) 
 \\
 &= & \CL_{R_\lambda}( \lambda(JY)) - (\CL_{R_\lambda}\lambda)(JY) \\
& = &  0-0 = 0
\eeastar
where we use $\lambda(J(\cdot)) = 0$ repeatedly for the first and the third equalities,
and then use the vanishing $\CL_{R_\lambda} \lambda = 0$ for the third. Combining the
above calculations, we obtain
\be\label{eq:symmetry-check}
\langle (\CL_{R_\lambda}J)Y,Z \rangle = d\lambda((\CL_{R_\lambda}J)Y, JZ)
\ee
Finally repeatedly applying the Leibniz rule to the Lie derivative $\CL_{R_\lambda}$
and the identity $d\lambda(J(\cdot), J(\cdot)) = d\lambda(\cdot, \cdot)$,
we rewrite the right hand side of \eqref{eq:symmetry-check} into
\beastar
&{}&  d\lambda((\CL_{R_\lambda}J)Y, JZ) \\
& = &  d\lambda( \CL_{R_\lambda}(JY) -   J(\CL_{R_\lambda}Y), JZ) \\
& = &  d\lambda(\CL_{R_\lambda}(JY),JZ)  -  d\lambda( J(\CL_{R_\lambda}Y), JZ) \\
 & = &  \big(\CL_{R_\lambda}(d\lambda( JY, JZ)) -   d\lambda (JY,  \CL_{R_\lambda}(JZ))\big)   
 - d\lambda\left(J\CL_{R_\lambda}Y, JZ\right)\\
& = &  \big( \CL_{R_\lambda}(d\lambda( JY, JZ)) -   d\lambda (JY,  (\CL_{R_\lambda}J)Z))   
 -   d\lambda (JY,  J \CL_{R_\lambda}Z) \big) 
  - d\lambda\left(J\CL_{R_\lambda}Y, JZ\right) \\
 & = &  \CL_{R_\lambda}(d\lambda(Y, Z)) - d\lambda(JY,  (\CL_{R_\lambda}J) Z)  
 - d\lambda\left(Y, \CL_{R_\lambda} Z\right) 
 -  d\lambda\left(\CL_{R_\lambda}Y, Z\right) \\
& = & \big( \CL_{R_\lambda}(d\lambda(Y, Z))   
 - d\lambda\left(Y, \CL_{R_\lambda} Z\right) 
 -  d\lambda\left(\CL_{R_\lambda}Y, Z\right) \big) - d\lambda(JY,  (\CL_{R_\lambda}J) Z)\\
& = &  \CL_{R_\lambda}(d\lambda)(Y, Z) - d\lambda(JY,  (\CL_{R_\lambda}J) Z)  \\
& =  & - d\lambda(JY,  (\CL_{R_\lambda}J) Z) 
=  d\lambda(Y,  J(\CL_{R_\lambda}J) Z) \\
& = & \langle Y,  (\CL_{R_\lambda}J) Z \rangle - \lambda(Y) \cdot \lambda(J(\CL_{R_\lambda}J)Z) 
=  \langle Y,  (\CL_{R_\lambda}J) Z \rangle.
 \eeastar
Combining this with \eqref{eq:symmetry-check}, we  have proved \eqref{eq:symmetry}.
This finishes the proof.
\end{proof}

\begin{rem}\label{rem:CLRJ=0} Using the vanishing
$\CL_{R_\lambda}\lambda=0$ and $ \CL_{R_\lambda}d\lambda= 0$, it is not hard to see that,
the Reeb vector field is a Killing vector field with respect to the triad metric if and only if $J$ satisfies $\CL_{R_\lambda}J=0$,
i.e., the tuple $(J,R_\lambda,\lambda, g)$ is K-contact in the sense of \cite{blair}.
(See Proposition \ref{prop:nablaYRlambda} below for further relevant comment.)
In general, this is a strong additional requirement. For example, for $3$-dimensional contact manifolds, 
it is equivalent to the Sasakian condition.  
\end{rem}

\subsection{Definition of contact triad connection}

To define the contact analog to the canonical connection of the
case of almost K\"ahler manifolds, we note that 
for any contact manifold $(M,\xi)$, the leaf space of the Reeb foliations of
contact form $\lambda$ becomes a naturally an (non-Hausdorff) almost
K\"ahler manifold $(\widehat M, \widehat{d\lambda}, \widehat J_\xi)$.
This picture motivates the following definition of
canonical connection of the contact triads $(M,\lambda,J)$.
We recall
$$
TM =  \xi \oplus \R\langle R_\lambda \rangle.
$$
We denote by $\pi = \pi_\lambda: TM \to \xi$ the associated projection. 

We reorganize the axioms of contact triad connection in different order 
from that of Theorem \ref{thm:mainI} as follows.

\begin{defn}[Contact triad connection]\label{defn:triad-connection} 
Let $(M,\lambda,J)$ be any contact triad of
contact manifold $(M,\xi)$. Let $g$ be the associated triad metric.
Then there exists a unique affine connection $\nabla$ that has the
following properties:
\begin{enumerate}
\item $\nabla$ is a Riemannian connection of the triad metric.
\item $\nabla_{R_\lambda} R_\lambda = 0$ and $\nabla_Y R_\lambda\in \xi$, for $Y\in \xi$.
\item For $Y\in \xi$, we have the following
$$
\del^\nabla_Y R_\lambda:= \frac12(\nabla_Y R_\lambda- J\nabla_{JY} R_\lambda)=0.
$$
\item The torsion tensor $T$ of $\nabla$ satisfies $T(R_\lambda, Y)=0$ for all $Y \in TQ$.
\item $\nabla^\pi := \pi \nabla|_\xi$ defines a Hermitian connection of the vector bundle
$\xi \to Q$ with Hermitian structure $(d\lambda, J)$.
\item The $\xi$ projection, denoted by $T^\pi: = \pi T$, of the torsion $T$
is of $(0,2)$-type in its complexification, i.e., satisfies the property \eqref{eq:TJYYxi}.
$$
T^\pi(JY,Y) = 0
$$
for all $Y$ tangent to $\xi$.
\end{enumerate}
We call $\nabla$ the contact triad connection.
\end{defn}

The following existence theorem is proved by Wang and the present author
in \cite{oh-wang:connection}. We will later provide a much more canonical and conceptual
proof than the original proof given therein by adapting Kobayashi's  moving frame method 
used in his proof of the unique existence theorem of a canonical connection
on almost Hermitian manifold \cite{kobayashi:canonical}.

\begin{thm}[\cite{oh-wang:connection}]\label{thm:connection-existence} There exists a unique connection $\nabla$ for $(M,\lambda,J)$
that satisfies the axioms given in Definition \ref{defn:triad-connection}.
\end{thm}

\begin{rem} More generally, it is proved in \cite{oh-wang:connection} that 
for any constant $c \in \R$, the same unique existence theorem holds by
replacing Condition (3) of Definition \ref{defn:triad-connection}  by
\be\label{eq:(3)-c}
\nabla_{JY} R_\lambda+J\nabla_Y R_\lambda = c \,Y
\ee
for $Y \in \xi$.
\end{rem}

\subsection{Properties of contact triad connection}

In this subsection, we collect various identities that contact triad connections satisfy.
All these identities turn out to dramatically simplify the a priori rather complex tensorial calculations
entering in the a priori estimates \cite{oh-wang:CR-map1,oh-wang:CR-map2,oh-yso:index}, 
the derivations of the discernable formulae for the linearization operator \cite{oh:contacton}
and the asymptotic operator \cite{kim-oh:asymp-analysis} of contact instantons and others.
Without the usage of the contact triad connection, all these tensorial calculations would
result in some formulae too complicated to discern, at least in Wang and the present
author's study of geometric analysis of contact instantons in \cite{oh-wang:CR-map1,oh-wang:CR-map2}.

First, it follows from definition that the contact triad connection $\nabla$ canonically induces
a Hermitian connection $\nabla^\pi$ for the Hermitian vector bundle $(\xi, J, g_\xi)$ over $M$.

Moreover, the following functorial properties of the contact triad connection were
proved in \cite{oh-wang:connection} and play an important role in the derivation of explicit
formulae of the linearized operator, its Fredholm analysis \cite{oh:contacton} and of the asymptotic operator 
\cite{kim-oh:asymp-analysis}, and their spectral analysis. 

\begin{prop}[Naturality] Let $(M,\lambda,J)$ be a contact triad, and 
 $\nabla$ be the contact triad connection thereof. Then 
 \begin{enumerate}
\item  For any diffeomorphism $\phi: Q \to Q$,
the pull-back connection $\phi^*\nabla$ is the triad connection
associated to the triad $ (M,\phi^*\lambda,\phi^*J)$ associated to
the pull-back contact structure $\phi^*\xi$.
\item
In particular if $\phi$ is contact, i.e., $d\phi(\xi) \subset \xi$, then
$ (M, \phi^*\lambda, \phi^*J)$ is a contact triad of $\xi$ and $\phi^*\nabla$
the contact triad connection $  (M,\xi)$.
\end{enumerate}
\end{prop}
The axiom that 
$$
\nabla_X R_\lambda \in \xi
$$
for all $X \in TM$ plays a fundamental role in the aforementioned
a priori estimates of the $\pi$-energy for both closed and open string contexts.
(See \cite{oh-wang:CR-map1,oh-wang:CR-map2} and \cite{oh-yso:index}, respectively.)

\begin{rem}
\begin{enumerate}
\item We also mention that the property $\nabla_{R_\lambda}J=0$ holds  
which is nothing but the vanishing
$$
\nabla_{R_\lambda}(JY)-J\nabla_{R_\lambda}Y=0
$$
for all $Y \in \xi$ which is embedded into the combination of the vanishing $J R_\lambda = 0$ 
and the Axiom of  the Hermitian property of the triple $(\xi,g|_\xi,J_\xi)$ of the complex vector bundle 
$(\xi, J_\xi) \to M$.
\item 
Similar vanishing $\nabla_{R_\lambda}^{LC} J = 0$ is known to
also hold for the Levi-Civita connection \cite[Corollary 6.1]{blair}. This can be easily derived from 
the vanishing $\nabla_{R_\lambda} J = 0$ for the triad connection $\nabla$ 
and the following identity 
$$
\nabla_{R_\lambda}Y = \nabla_{R_\lambda}^{LC}Y - \frac12 JY
$$
which is proved in \cite[Lemma 4.6]{kim-oh:asymp-analysis}.
(See  \cite[Section 4.2]{kim-oh:asymp-analysis} for detailed discussion on this property.)
\end{enumerate}
\end{rem}

The following identities are also very useful to perform tensorial
calculations in the study of a priori elliptic estimates and
in the derivation of the linearization formula and asymptotic analysis of
contact instantons \cite{oh-wang:CR-map1,oh-wang:CR-map2}, \cite{oh-yso:index}
and \cite{kim-oh:asymp-analysis}.

\begin{prop}\label{prop:nablaYRlambda}
Let $\nabla$ be the contact triad connection. Then
 for any vector field $Y$ on $Q$,
\be\label{eq:nablaRlambda}
\nabla_Y R_\lambda = \frac{1}{2}(\CL_{R_\lambda}J)JY.
\ee
In particular the linear operator $Y \mapsto \nabla_Y R_\lambda$ is a 
symmetric operator with respect to the triad metric.
\end{prop}
\begin{proof} The statement of the symmetry of the operator follows from 
Proposition \ref{prop:symmetry}.  Therefore we will focus on the 
proof of the formula \eqref{eq:nablaRlambda}.

First, we recall the Hermitian connection property,  and torsion property
$T^\pi|_\xi$ has vanishing $(1,1)$ part.

 Combining $\nabla_{R_\lambda}J = 0$ with the torsion property, we compute 
\beastar
0 & = & T(R_\lambda, JY) = \nabla_{R_\lambda} (JY) - \nabla_{JY}(R_\lambda) - [R_\lambda, JY] \\
& = &  J \nabla_{R_\lambda} Y - \nabla_{JY}(R_\lambda) - J [R_\lambda, Y] 
- (\CL_{R_\lambda}J)(Y)  \\
& = &  J (\nabla_{R_\lambda} Y - \nabla_Y R_\lambda -  [R_\lambda, Y] ) 
 + (J \nabla_Y R_\lambda - \nabla_{JY}(R_\lambda)) - (\CL_{R_\lambda}J)(Y) \\
 & = & J T(R_\lambda, Y) + 2J \nabla_{Y} R_\lambda  - (\CL_{R_\lambda}J) (Y) \\
 & = & 2J \nabla_{Y} R_\lambda  - (\CL_{R_\lambda}J) (Y)
 \eeastar
where we use the vanishing $\nabla_{R_\lambda}J = 0$ for the second equality,
 rearranging for the third and Axiom (3) and (4) for the fourth and the fifth in turn.
Therefore we obtain
$$
\nabla_{R_\lambda} Y =  - \frac12 J (\CL_{R_\lambda}  J) Y = \frac12  (\CL_{R_\lambda}  J) J Y
$$
where we use $J \CL_{R_\lambda}J = - (\CL_{R_\lambda}J)J$.
\end{proof}

\begin{rem}
We can relax Condition (3) by replacing it to
$$
\nabla_{JY}R_\lambda+J\nabla_Y R_\lambda = cY,
$$
for any given real number $c$. This way we shall have one-parameter family of
affine connections parameterized by the real number $c$ each of which satisfies Condition (1)-(6) except (3)
replaced by \eqref{eq:(3)-c}.

With $c$ fixed, i.e., under Condition \eqref{eq:(3)-c},
$\nabla_Y R_\lambda$ satisfies
$$
\nabla_Y R_\lambda=-\frac{1}{2}cJY+\frac{1}{2}(\CL_{R_\lambda}J)JY.
$$
\end{rem}

We also state the following useful property for the completeness' sake although 
it is not explicitly used in the current paper. It is a consequence of Axiom (2) which 
is useful for the analysis of contact instantons \cite{oh-wang:CR-map1,oh-yso:index}. 

\begin{prop}\label{prop:nablaRlambda=0} We have $\nabla_{R_\lambda} \lambda = 0$ and
\beastar
 \nabla_{R_\lambda}(d\lambda) 
& = & \left\langle \CL_{R_\lambda}J(\cdot), \cdot \right \rangle \\
& = &d\lambda((\CL_{R_\lambda}J)(\cdot),J(\cdot)).
\eeastar
\end{prop}
\begin{proof} For the fist vanishing, we evaluate $(\nabla_{R_\lambda} \lambda)(X)$.
If $X= R_\lambda$, we derive
$$
(\nabla_{R_\lambda} \lambda)(R_\lambda) = \nabla_{R_\lambda}(\lambda(R_\lambda)) 
- \lambda(\nabla_{R_\lambda}R_\lambda) = 0 - 0 = 0.
$$
For $Y \in \xi$, i.e., $Y$ is a vector field tangent to $\xi$, we evaluate
$$
(\nabla_{R_\lambda} \lambda)(Y) = \nabla_{R_\lambda}(\lambda(Y)) 
- \lambda(\nabla_{R_\lambda}Y) = - \lambda(\nabla_{R_\lambda}Y).
$$
Then using the vanishing $T(R_\lambda,Y) = 0$ and \eqref{eq:nablaRlambda}, we rewrite the 
argument of the last term into
$$
\nabla_{R_\lambda}Y = \nabla_Y R_\lambda + [R_\lambda,Y] 
= \frac12 (\CL_{R_\lambda} J) JY + \CL_{R_\lambda}Y 
$$
which is contained in $\xi$ since $Y \in \xi$ by assumption. This proves $ \lambda(\nabla_{R_\lambda}Y) = 0$
for all $Y \in \xi$.
Combining the two vanishing, we obtain $\nabla_{R_\lambda}\lambda = 0$.

For the second vanishing, we first rewrite
\bea\label{eq:nabladlambda}
&{}& \nabla_{R_\lambda}(d\lambda)(X_1,X_2) \nonumber \\
 & = & \nabla_{R_\lambda}(d\lambda(X_1, X_2))
- d\lambda(\nabla_{R_\lambda} X_1, X_2) -  d\lambda(X_1, \nabla_{R_\lambda}X_2) \nonumber\\
& = &  \nabla_{R_\lambda}(d\lambda(X_1, X_2))
- \big( d\lambda(\nabla_{R_\lambda} X_1, X_2) + d\lambda(X_1, \nabla_{R_\lambda}X_2)\big).
\eea
When $X_1 = R_\lambda$ and $X_2 = Y \in \xi$ or vice versa,  the right hand side vanishes. 
Therefore it remains to consider the case with
$X_1 =: Y, \, X_2 = : Z \in \xi$.  In that case, we see
\be\label{eq:first}
\nabla_{R_\lambda}(d\lambda(Y, Z)) = \CL_{R_\lambda}(d\lambda(Y,Z))
\ee
for the first term. For the rest of the right hand side of
\eqref{eq:nabladlambda}, using $R_\lambda \intprod T = 0$, we evaluate
 the sum inside the outside parenthesis into
\bea\label{eq:2nd+3rd}
&{}& d\lambda(\nabla_{R_\lambda} Y, Z) +  d\lambda(Y, \nabla_{R_\lambda} Z)  \nonumber \\
& = & d\lambda(\nabla_Y R_\lambda + [R_\lambda, Y], Z)  
+ d\lambda(Y, \nabla_Z R_\lambda + [R_\lambda, Z]) \nonumber\\
& = & d\lambda(\nabla_Y R_\lambda ,Z) 
+ d\lambda(Y, \nabla_Z R_\lambda) + \big(d\lambda ([R_\lambda, Y],Z)+ d\lambda(Y,[R_\lambda, Z])\big)
\nonumber\\
 & = & d\lambda\left(\frac12 \CL_{R_\lambda}J JY, Z\right) 
 + d\lambda\left(Y, \frac12 \CL_{R_\lambda}J JZ\right) + 
 \CL_{R_\lambda}(d\lambda(Y,Z)) \nonumber \\
& = & - \left\langle \frac12( \CL_{R_\lambda}J) JY, JZ\right\rangle 
 - \left \langle Y, \frac12 J ( \CL_{R_\lambda}J JZ)\right \rangle + 
 \CL_{R_\lambda}(d\lambda(Y,Z)) \nonumber \\
 & = &\frac12 \big(- \left\langle ( \CL_{R_\lambda}J) JY, JZ\right\rangle 
 + \left \langle JY, ( \CL_{R_\lambda}J) JZ\right \rangle\big) + 
 \CL_{R_\lambda}(d\lambda(Y,Z)) \nonumber \\
  & =  & -  \left \langle ( \CL_{R_\lambda}J) JY, JZ\right\rangle + \CL_{R_\lambda}(d\lambda(Y,Z))
  \nonumber\\
  & = &  \left \langle ( \CL_{R_\lambda}J) Y, Z\right\rangle + \CL_{R_\lambda}(d\lambda(Y,Z))
    \eea
where we use Proposition \ref{prop:nablaYRlambda},
the symmetry of $ \CL_{R_\lambda}J $ for the penultimate equality, and
$\CL_{R_\lambda}J J = - J \CL_{R_\lambda}J$ and $J$ being an isometry for the last equality.
Substituting \eqref{eq:first} and \eqref{eq:2nd+3rd} into \eqref{eq:nabladlambda}, we conclude
$$
\nabla_{R_\lambda}(d\lambda)(Y,Z) =   \left \langle ( \CL_{R_\lambda}J) Y, Z\right\rangle 
= d\lambda((\CL_{R_\lambda}J)Y,JZ)
$$
for all $Y, \, Z \in \xi$.  Combining the above calculations, we have finished the  proof.
\end{proof}

\begin{rem} 
We would like to alert readers and compare this lemma
with the vanishing $\CL_{R_\lambda}\lambda = 0, \, \CL_{R_\lambda}d\lambda = 0$ of Lie
derivatives which is an easy functorial consequence  of Cartan's magic formula and the defining 
properties of the Reeb vector field.
\end{rem}

Next we prove an interesting and useful torsion property of the triad connection.

\begin{prop} \label{prop:lambda(T)}
For the Reeb component of the torsion $T$, we have
\be\label{eq:lambdaT}
\lambda (T ) = d\lambda.
\ee
\end{prop}
\begin{proof} We consider two cases separately, one for $(X_1, X_2) = (R_\lambda, Y)$ with $Y \in \xi$
and the other for $X_1 = Y, \, X_2 = Z \in \xi$.
For the first case, we know $d\lambda(R_\lambda,Y) = 0$ and $\lambda(T(R_\lambda,Y)) = 0$
thanks to the vanishing $R_\lambda \intprod T = 0$.

For the second case, we evaluate
\bea\label{eq:dlambdaYZ}
d\lambda(Y, Z) & = & Y [\lambda(Z)] - Z(\lambda(Y)) - \lambda([Y,Z])  \nonumber\\
& = & 0 - 0 -\lambda([Y,Z]) = -\lambda([Y,Z]).
\eea
On the other hand, we have
\bea\label{eq:lambdaTYZ}
\lambda(T(Y, Z)) & = & \lambda(\nabla_Y Z - \nabla_Z Y - [Y,Z]) \nonumber\\
& = & \lambda(\nabla_Y Z) - \lambda(\nabla_{Z} Y) - \lambda([Y,Z]) \nonumber\\
& = & - (\nabla_Y \lambda)(Z) +  (\nabla_Z \lambda)(Y)   - \lambda([Y,Z])
\eea
by the definition of $T$ for the first equality and thanks to $\ker \lambda = \xi$ 
and the Lebnitz rule for the covariant derivative for the last two
equalities.

For the first two terms, we rewrite
\beastar
(\nabla_Y \lambda)(Z) & = & - \lambda(\nabla_Y Z) = - \langle R_\lambda, \nabla_Y Z \rangle 
=  \langle \nabla_Y R_\lambda, Z \rangle\\
& = & \frac12 \langle \CL_{R_\lambda}JJ Y,Z \rangle
= \frac12 \langle Y,   \CL_{R_\lambda}JJ Z \rangle = \nabla_Z\lambda(Y)
\eeastar
where we use the vanishing $\lambda|_\xi = 0$, Proposition  \ref{prop:nablaYRlambda}
and the symmetry of $\CL_{R_\lambda}JJ$.
Substituting the above into \eqref{eq:lambdaTYZ}, we have derived
$$
\lambda(T(Y,Z)) =   - \lambda([Y,Z]).
$$
Comparing this with \eqref{eq:dlambdaYZ}, we prove $d\lambda(Y,Z) = \lambda(T(Y,Z))$
for $Y, \, Z \in \xi$. Combining the above,
 we have finished the proof of Proposition \ref{prop:lambda(T)}.
\end{proof}

\section{Connection in almost contact frame and its first structure equations}
\label{sec:in-complexframe}

We first reformulate the defining conditions of the triad connection in terms of
the moving Darboux frame.

\subsection{Structure equations in the Darboux frame}

Choose a local orthogonal Darboux  frame of $TM = \xi \oplus\R \langle R_\lambda \rangle $ given by
$$
\CD: = \{ E_1,\cdots, E_n, F_1,\cdots, F_n,  R_\lambda\}, \quad F_i: = J E_i
$$
and denote its dual co-frame by
$$
e^1,\cdots, e^n, f^1,\cdots, f^n, \, \lambda.
$$
Assume the connection matrix of one-forms is $(\Gamma^i_j)$, $i, j=0, 1,..., 2n$
with
$$
\Gamma^i_k : = \sum_k \Gamma^i_{kj} e^j + \sum_k \Gamma^i_{k, n+j} f^j
$$
and we write the first structure equations as follows.
\be\label{eq:structure-equation-real}
\begin{cases}
d\lambda = -\Gamma^0_0\wedge \lambda- \sum_k \Gamma^0_k\wedge e^k- \sum_k\Gamma^0_{n+k}\wedge f^k+T^0\\
d e^j = -\Gamma^j_0\wedge \lambda- \sum_k \Gamma^j_k\wedge e^k- \sum_k \Gamma^j_{n+k}\wedge f^k+T^j\\
d f^j = -\Gamma^{n+j}_0\wedge \lambda- \sum_k \Gamma^{n+j}_k\wedge e^k-\sum_k \Gamma^{n+j}_{n+k}\wedge f^k+T^{n+j}
\end{cases}
\ee
Throughout the section, if not stated otherwise, we let  $i$, $j$ and $k$ take values from $1$ to $n$.
A straightforward calculation leads to
the following vanishing results of the Christoffel symbols 
which are derived in \cite{oh-wang:connection} and its arXiv version.

\begin{lem}[p.10 \cite{oh-wang:connection}]\label{lem:Gamma00=0} Let $\nabla$ be any contact triad $c$-connection. Then 
$\Gamma^0_{00} =0$ and 
$$
\Gamma^a_{00} = 0 = \Gamma^0_{b0}\quad \forall a,\, b = 1, \cdots 2n
$$
for any orthogonal Darboux frame. In particular we have $\Gamma^0_0 = 0$.
\end{lem}
\begin{proof} These are immediate consequences of the properties $\nabla_{R_\lambda}R_\lambda =0$
and Axiom (2) of Definition \ref{defn:triad-connection}.
\end{proof}

The following is also proved in \cite{oh-wang:connection} (or its arXiv version).

\begin{lem}[Equations (15) \& (16) \cite{oh-wang:connection}] \label{lem:alpha0-antiholomorphic-c}
Axiom (3) is equivalent to
$$
\Gamma^k_{j0} + \Gamma^{n+k}_{n+j,0}=0 , \quad \Gamma^{n+k}_{j,0} - \Gamma^k_{n+j,0} = - c\, \delta_{j,k}
$$
for all $j, \, k = 1, \cdots n$
in the real structure equations \eqref{eq:structure-equation-real}.
\end{lem}
 
We will analyze each condition in Definition \ref{defn:triad-connection} and show how they set down the matrix of connection one-forms.

\subsection{Structure equation in the almost-contact frame}

We now introduce a complex-valued one-forms 
\be\label{eq:etai}
\eta_i : = \frac1{\sqrt{2}}(E_i- \sqrt{-1}  F_i), \quad \overline \eta_i : = \frac1{\sqrt{2}}(E_i+ \sqrt{-1}  F_i)
\ee
for $i = 1, \ldots, n$.

\begin{defn}[Almost contact frame] Let $(M,\lambda, J)$ be a contact triad equipped with triad metric.
 We call the complex frame 
\be\label{eq:ac-frame}
\{\eta_1,\cdots, \eta_n, \overline \eta_1, \cdots, \overline \eta_n\, R_\lambda\} = :\CD_\C
\ee
the \emph{almost contact frame} associated to the Darboux frame $\CD$, and just call any
such frame an almost contact frame when the Darboux frame $\CD$ is not specified.
\end{defn}
Complexification of $TM$ and complex linear extension $g_\C = g \otimes \C$
of the contact triad metric  $g$ also provide  the complex vector bundle
$T^*M \otimes \C$ with the associated dual frame 
$$ 
\{\theta^1, \cdots, \theta^n, \overline \theta^1, \cdots, \overline \theta^ n\} \oplus \C \langle \lambda \rangle.
$$
We will also call it the  \emph{almost contact frame} of the contact triad.
We  have the decomposition
\be\label{eq:dual-decomposition}
T^*M \otimes \C = (\xi^* \otimes \C) \oplus \C \langle \lambda \rangle = 
 \xi^{*(1,0)} \oplus \xi^{*(0,1)} \oplus \C \langle \lambda \rangle.
\ee

We now express a connection on the complex subbundle 
$$
\xi^{(1,0)} \oplus \C \langle R_\lambda \rangle \subset TM \otimes \C
$$
over $M$  in terms of the frame  \eqref{eq:ac-frame}, i.e., in terms of
$$
\C\{\eta_1,\cdots, \eta_n\} \oplus \C \{R_\lambda\}, \quad \eta_i : = \frac1{\sqrt{2}}(E_i- \sqrt{-1}  F_i)
$$
where $\xi^{(1,0)}$ is the $(1,0)$ part of $\xi_\C = (\xi,J) \otimes \C$.  Then $g_\C$ induces a Hermitian
metric on the complex vector bundle 
$\xi^{(1,0)} \oplus \{0\} \subset \xi^{(1,0)} \oplus \C\langle R_\lambda \rangle$
where $\xi^{(1,0)}$ is the $(1,0)$ part of $\xi_\C = \xi \otimes \C$. 

We have explicit expressions of the metric $g_\C$ on $\xi^{(1,0)} \oplus \C\langle R_\lambda \rangle$
as the bilinear form 
$$
g = \sum_k \theta^k  \otimes\overline \theta^k + \lambda \otimes \lambda
$$
which restricts to a Hermitian metric on $\xi^{(1,0)}$
where the associated Hermitian frame is expressed as 
$$
\theta^i = \frac1{\sqrt{2}} (e^i + \sqrt{-1} f^i), \, \overline \theta^i 
= \frac1{\sqrt{2}} (e^i - \sqrt{-1} f^i).
$$
Furthermore we can express the 2-form
\be\label{eq:dlambda}
d\lambda = g(J \cdot, \cdot) = \sqrt{-1}\sum_k \theta^k \wedge \overline \theta^k
\ee
in terms of the  almost contact frame
\be\label{eq:complex-frame}
\C\{\eta_1,\cdots, \eta_n\} \oplus \C \{R_\lambda\}, \quad \eta_i : = \frac1{\sqrt{2}}(E_i- \sqrt{-1}  F_i).
\ee

Any connection on the complex vector bundle $\xi^{(1,0)} \oplus \C \langle R_\lambda \rangle \to M$
can be expressed in terms of the frame \eqref{eq:complex-frame}.
We denote
its dual frame by
$$
\{\theta^1, \cdots, \theta^n\} \oplus \{\lambda\}
$$
where we recall that $\theta^i$ is 1-forms of $(1,0)$-type and $\lambda$ is a real one-form.
The connection form of the complex vector bundle $\xi^{(1,0)} \oplus \C \langle R_\lambda \rangle \to M$
associated to the above frame is a matrix-valued one-form  consisting of a skew-Hermitian matrix
$$
\omega_j^i, \quad i, \, j = 1,\cdots, n
$$
and of complex vectors
$$
\alpha^0 = (\alpha^0_j) \quad \beta_0^t = (\beta_0^i), \quad i,\, j = 0, \cdots, n
$$
of one-forms.

\begin{notation} We recall the splitting $TM = \xi \oplus \R \langle R_\lambda \rangle$.
With respect to this splitting, we have the associated splitting
$$
T^*M = \R \langle R_\lambda \rangle^\perp \oplus \xi^\perp \cong \xi^* \oplus \R \langle \lambda \rangle.
$$
For a given one-form $\gamma$, we denote by $\gamma^\pi$ its $\pi$-component  given by
\be\label{eq:alpha-pi}
\gamma^\pi = \gamma- \gamma(R_\lambda) \lambda.
\ee
 By definition, we have $\gamma = \gamma^\pi + \gamma(R_\lambda) \lambda$. We call $\gamma^\pi$
and $\gamma^\perp := \gamma(R_\lambda) \lambda$ the \emph{horizontal} and the \emph{vertical} 
or the Reeb component of $\gamma$.
For example, we will have
$$
(e^i)^\pi = e^i, \,\,  (f^j)^\pi = f^j
$$
and their vertical components being zero,
 and $\lambda^\perp = \lambda$ and its horizontal component being zero,
when we apply the practice to the basis elements of the almost contact frame.
\end{notation}

In terms of the given almost contact frame $\CD_\C$,
the full complex-valued connection  can be expressed as a matrix-valued one-form
denoted by $\Omega = (\Omega_j^i)$ given by
\be\label{eq:Omega}
\Omega = \left(\begin{matrix} \omega & \alpha^{0\pi} \\
(\beta^{0\pi})^t & 0 \end{matrix}\right)
\ee
and its complex torsion  $\Theta := \Pi' T_\C$$\Theta$ as the matrix-valued two form, which is denoted by
$$
\Theta = (\Theta^1, \cdots, \Theta^n, \Theta^0)
$$
in matrix, where its components $\Theta^i$ are given by
\be\label{eq:Theta-defn}
\begin{cases}
\Theta^0 := \lambda(T_\C) = \lambda(T_\C(\cdot, \cdot))\\
 \Theta^i : = \theta^i(T_\C) = \theta^i(T_\C(\cdot,\cdot)) \quad \text{\rm for } \ i=1, \cdots, n
 \end{cases}
\ee
by definition, which we regard ordinary complex-valued two-forms.

Then the first structure equations of a connection in the given almost contact frame 
is given by 
\bea
\Theta^i & = &d\theta^i+  \sum_{k=1}^n \omega_k^i \wedge \theta^k + \beta^i_0 \wedge \lambda, \quad i = 1, \ldots, n,
\label{eq:torsion-Thetai} \\
\Theta^0 & = & d\lambda + \sum_{k=1}^n \alpha_k^0 \wedge \theta^k + \alpha_0^0 \wedge \lambda.
\label{eq:torsion-Theta0}
\eea
We will often just write
$$
\beta^i_0 = \beta^i
$$
omitting the unnecessary subindex 0 in the later calculations.

\section{Contact triad connection in the almost contact frame}
\label{sec:connection-in-acframe}

Starting from this section and the rest, we will study contact triad connection, more generally 
 with $c$ not necessarily being 0.
Not to confuse readers, we will call the corresponding connection \emph{$c$-contact triad connection}, 
or simply \emph{$c$-connection}. We first mention the following immediate translations of Axiom (2) of Definition \ref{defn:triad-connection} and the skew-Hermitian property of $\Omega$.

\begin{lem}\label{lem:Axiom(2)(4)} 
Axiom (2) of Definition \ref{defn:triad-connection} and the skew-Hermitian property of $\Omega$ 
as above are equivalent to 
\be\label{eq:alpha00}
\alpha_0^0 = 0
\ee
and 
\be\label{eq:alphabar}
 \beta^k  = - \overline{\alpha}^0_k \quad k = 1, \cdots n
\ee
respectively.
\end{lem}
We also have the following translation of Axiom (3), or more generally \eqref{eq:(3)-c}.
Its proof displays an interesting interaction between the axiom of contact triad connection
and the structure equations thereof
\emph{in terms of the almost contact frame}. 

\begin{prop}\label{prop:axiom-3} The equation \eqref{eq:(3)-c} is equivalent to
\be\label{eq:alpha0-antiholomorphic-c}
J\cdot \alpha^0_k = - \sqrt{-1} \alpha^{0}_k + c\, \theta^k, \quad k = 1, \cdots, n.
\ee
In particular, when $c = 0$, $\alpha^0_k$ is anti-holomorphic.
\end{prop}
\begin{proof} First, by definition of the torsion $T$ and \eqref{eq:lambdaT}
in turn, we derive
$$
\Theta^0 = \lambda(T) = d\lambda.
$$
Therefore this combined with \eqref{eq:torsion-Theta0} and $\alpha_0^0 = 0$ gives rise to
\be\label{eq:Theta0=dlambda}
\sum_{k=1}^n \alpha_k^{0} \wedge \theta^k = \Theta^0 - d\lambda = 0.
\ee
We decompose
$$
\alpha_k^0 =  (\alpha_k^0)^\pi + \alpha_k^{0\perp}  =  (\alpha_k^{0\pi})^{(1,0)} +  (\alpha_k^{0\pi})^{(0,1)}
+  \alpha_k^{0\perp}
$$
for $k = 0, \ldots, n$. 
Writing
$$
 \alpha_k^{0\perp} = \alpha_k^0(R_\lambda) \lambda
$$
and applying Axiom (2) (or Lemma \ref{lem:Gamma00=0} directly), we derive
\be\label{eq:alphak00=0}
 \alpha_k^{0\perp} = 0
 \ee
hence $\alpha^0_k = \alpha^{0\pi}_k$.
Substitution of these into \eqref{eq:Theta0=dlambda}
and a simple type analysis gives rise to
\bea
\sum_{k=1}^n (\alpha_k^{0\pi})^{(1,0)} \wedge \theta^k & =  & 0, \label{eq:(1,0)thetak}\\
\sum_{k=1}^n (\alpha_k^{0\pi})^{(0,1)} \wedge \theta^k & = & 0 \label{eq:(0,1)thetak}\\
\eea
We note that the set $\{\theta^k\}_{k=1}^n$ (resp. $\{\bar \theta^k\}_{k=1}^n\}$) is 
linearly independent on the $\C$-vector space $\xi^{(1,0)}$ (resp. on $\C$-vector space $\xi^{(0,1)}$).

Applying Cartan's lemma to \eqref{eq:(1,0)thetak} (see \cite[Problems 11, Chapter 7]{spivakI},
for example) to the pair of sets $\{\theta^k\}_{k=1}^n$ and
$$
 \left\{(\alpha_k^{0\pi})^{(1,0)} \right\}_{k=1}^n 
 $$
 respectively, we obtain 
\be\label{eq:alphak0}
(\alpha_k^{0\pi})^{(1,0)} = \sum_{i =1} b_{ki} \theta^i 
\ee
for some $b_{ki}$ satisfying $b_{ik} = b_{ki}$.  

On the other hand, noting $\overline \theta_i \wedge \theta_j$ are
linearly independent and substituting the formula
$$
\left\{(\alpha_k^{0\pi})^{(0,1)}\right \}_{k=1}^n
$$
into \eqref{eq:(0,1)thetak}, we obtain $c_{k\bar{i}} = 0$ for all $k,\, i$. Therefore we have obtained
$(\alpha_k^{0\pi})^{(0,1)}=0$. 
Combining this vanishing with  \eqref{eq:alphak0}, \eqref{eq:alpha00} we have shown that
$$
\alpha_k^0 = \alpha_k^{0\pi} =  (\alpha_k^{0\pi})^{(1,0)} = 
\sum_{ i =1} b_{ki} \theta^i.
$$
This reduces the analysis of
$\alpha_k^0$ to that restricted to 
$$
(\xi\otimes \C)^{(1,0)} \subset T^{(1,0)}M \subset TM \otimes \C
$$
which is now in order.

By definition of the dual one-forms $e^i, \, f^j$, of the connection forms $\alpha_k^0$, the vanishing 
$\alpha_0^0 = 0$
and the reality of $\nabla_Y R_\lambda$, we derive
\beastar
\nabla_Y R_\lambda & = & \sum_{k=1}^n \alpha^0_k(Y) \eta_k 
+ \sum_{k=1}^n \overline \alpha^0_k(Y)  \overline \eta_k\\
J\nabla_{JY} R_\lambda & = &  \sum_{k=1}^n \alpha^0_k(JY) \sqrt{-1} \eta_k
+ \sum_{k=1}^n \overline \alpha^0_k(JY) (-\sqrt{-1}) \overline \eta_k.
 \eeastar
 By subtracting the second from the first, we obtain
 \beastar
 \nabla_Y R_\lambda - J\nabla_{JY} R_\lambda & = &  \sum_{k=1}^n 
 \left(\alpha^0_k(Y) - \sqrt{-1} \alpha^0_k(JY)\right) \eta_k \\
&{}& + \sum_{k =1}^n \left(\overline \alpha^0_k (Y)  + \sqrt{-1} \overline \alpha^0_k(JY))\right) \overline \eta_k.
\eeastar
Therefore 
\bea\label{eq:2JdelnablaY}
&{}&
J(\nabla_Y R_\lambda - J\nabla_{JY} R_\lambda) \nonumber\\
& = &  \sum_{k=1}^n 
 \left(\alpha^0_k(Y) - \sqrt{-1} \alpha^0_k(JY)\right) \sqrt{-1} \eta_k \nonumber \\
&{}& + \sum_{k =1}^n \left(\overline \alpha^0_k(Y)  + \sqrt{-1} \overline \alpha^0_k(JY))\right) (-\sqrt{-1}) 
\overline \eta_k \nonumber \\
& = &  \sum_{k=1}^n 
 \left(\alpha^0_k(JY) +  \sqrt{-1}\alpha^0_k(Y) \right)  \eta_k
+ \sum_{k =1}^n \left( \overline \alpha^0_k (JY)) - \sqrt{-1}\overline \alpha^0_k(Y) \right) \overline \eta_k
\eea
Recall that we can express any \emph{real} vector as
$$
Y = \sum_k \theta^k(Y) \eta_k + \sum_k \overline \theta^k(Y) \overline \eta_k.
$$
Therefore evaluating $\theta_k$ against \eqref{eq:(3)-c}, which is equivalent to
$$
J(\nabla_Y R_\lambda - J\nabla_{JY} R_\lambda) = c\, Y,
$$
 we derive
\be \label{eq:alpha-theta}
\alpha^0_k(JY) +  \sqrt{-1}\alpha^0_k(Y)  = c\, \theta^k (Y)
 \ee
from \eqref{eq:2JdelnablaY} for all $Y \in \xi$. Therefore we obtain
$$
J\cdot \alpha^0_k = - \sqrt{-1} \alpha_k^{0} + c \,\theta^k.
$$
This finishes the proof.
 \end{proof}

\begin{cor} Let $c = 0$. Then $\alpha^0_k$ satisfies 
$J\cdot \alpha^0_k = - \sqrt{-1} \alpha_k^{0}$, i.e., 
each $\alpha^{0}_k$ is a 1-form of $(0,1)$-type, and hence
$$
(\alpha_k^{0\pi})^{(0,1)} = \sum_{i=1}^n c_{i\bar k} \, \bar \theta^k = 0
$$
for each $k = 1, \ldots, n$. 
\end{cor}

Now the above discussion enables us to translate the condition on the real torsion $T$
laid out in Definition \ref{defn:triad-connection} to those of 
the complex torsion $\Theta = \Pi' T_\C$ using the expression \eqref{eq:Theta-defn} of its component
as follows. We denote the $\xi$-component (or $\pi$-component) of $\Theta$ by
\be\label{eq:PiTheta}
\Theta^\pi: = \pi \Theta 
\ee
for the projection $\pi = \pi_\lambda : TM \to \xi$, and decompose 
$$
\Theta^\pi = \Theta^{\pi(2,0)} + \Theta^{\pi(1,1)} + \Theta^{\pi(0,2)}
$$
as a $(\xi \otimes \C)^{(0,1)} = T_\C^{\pi(1,0)}M$-valued two form on $M$.

 Then the following is immediate by definition and from the above discussion.

\begin{prop}\label{prop:in-movingframe} Let $\nabla$ be any connection $\nabla g = 0$.
Then the following equivalences hold for any constant $c$:
\begin{enumerate}
\item Axiom \eqref{eq:(3)-c} is equivalent to \eqref{eq:alpha0-antiholomorphic-c}.
\item Axioms (2), (4) together is  equivalent to $R_\lambda \intprod \Theta^0 = 0$.
\item Conditions (1) and (5) are equivalent to $\Theta^\pi$ is skew-Hermitian.
\item Condition (6) is equivalent to the vanishing 
$(\Theta^i)^{(1,1)} = 0$ of the $(1,1)$-component of $\Theta^i$ for $i =1,\cdots, n$.
\end{enumerate}
\end{prop}
This completes translation of the defining conditions given in Definition \ref{defn:triad-connection}
into the languages of the almost contact moving frame $\CD_\C$ given in \eqref{eq:ac-frame}.

\section{Construction of contact triad connection}
\label{sec:construction}

Choose a local orthogonal Darboux  frame of $TM = \xi \oplus\R \langle R_\lambda \rangle $ given by
$$
 E_1,\cdots, E_n, F_1,\cdots, F_n,  R_\lambda, \quad F_i: = J E_i
 $$
and denote its dual co-frame by
$$
e^1,\cdots, e^n, f^1,\cdots, f^n, \, \lambda.
$$
By Proposition \ref{prop:in-movingframe}, it is enough to construct a complex connection
$\Omega$ on $\xi^{(1,0)} \oplus \C \langle R_\lambda \rangle$ such that its torsion $\Theta$ satisfies
the required properties laid out in the previous section.

Let  $\{\theta^1, \cdots, \theta^n, \lambda\} $ be the almost contact frame 
of $(\xi^*)^{(1,0)} \oplus \C \langle R_\lambda \rangle$ as in the previous section.
The first structure equation of a \emph{general} connection   is given by
\be\label{eq:first-dlambda}
d\lambda =  - \sum_k \alpha_k^0 \wedge  \theta^k - \alpha_0^0 \wedge \lambda + \Theta^0
\ee
and
\be\label{eq:first-dtheta}
d\theta^i = - \sum_{j=1}^n \omega_j^i \wedge \theta^j - \sum_{i=1}^n \beta_0^i \wedge \lambda + \Theta^i
\ee
where $\theta$ and $\Theta$ denote column vectors of 1-forms $\lambda = \theta^0,\cdots, \theta^n$ and
2-forms $\Theta^0,\dots, \Theta^n$. 

\subsection{Uniqueness of contact triad connection}

We start with the uniqueness of the above constructed triad connection for each given constant $c$. 
Recall  that such a connection satisfies
$$
\alpha_0^0= 0.
$$
Then combining this with and that of $\Theta^j$ for $j = 1, \cdots, n$ which is \eqref{eq:first-dtheta},
we show that the first structure equation is reduced to 
\be\label{eq:Theta0i}
\begin{cases}
\Theta^i = d\theta^i + \sum_{j=1}^n \omega_j^i \wedge \theta^j + \beta^i \wedge \lambda\\
\Theta^0 = d\lambda  + \sum_k \alpha_k^0 \wedge  \theta^k + \alpha^0 \wedge \lambda.
\end{cases}
\ee
From now on, we assume that the first structure equations of the given $c$-connections have this form
in the given almost contact frame.

For any two connection forms $\omega = (\omega_i^j)_{1 \leq i, \, j \leq n}$ and $\widetilde \omega = (\widetilde \omega_i^j)$ 
of $c$-connection,
the above stated properties have the same $(0,2)$-component of $\Theta$ and $\widetilde \Theta$
because this component is precisely Nijenhuis integrability tensor of the leaf space $M/\sim$ of
Reeb trajectories, which is independent of the choice of connections. 

Let $\omega$ and $\widetilde \omega$ be the connection forms of $c$-connections so that 
they have vanishing $(1,1)$-components.
Then, it follows from the first equation of \eqref{eq:Theta0i} that
\be\label{eq:omega-tildeometa}
\sum_{j=1}(\omega_j^i - \widetilde \omega_j^i)\wedge \theta^j 
= (\Theta_j^i)^{(2,0)} - (\widetilde \Theta_j^i)^{(2,0)} = 0.
\ee
This implies $\omega_j^i - \widetilde \omega_j^i$ has bi-degree $(1,0)$. Since $\omega - \widetilde \omega$
 is skew-Hermitian, we have
 $$
 \omega - \widetilde \omega = - \overline{(\omega - \widetilde  \omega)^t}.
 $$
 Noting that the right hand side  is of bi-degree $(0,1)$, we conclude 
 \be\label{eq:omega-equal}
 \widetilde \omega = \omega.
 \ee
On the other hand, we have 
\be\label{eq:alpha0-equal}
\widetilde \alpha_0^0 = 0 = \alpha_0^0
\ee
 from $\nabla_{R_\lambda}R_\lambda = 0$.
Furthermore Proposition \ref{prop:axiom-3} implies that $\widetilde \alpha^0_k$ satisfies
$$
J\cdot \widetilde \alpha^0_k = - \sqrt{-1} \widetilde \alpha^{0\pi}_k + c\, \theta^k
$$
Therefore the difference $\widetilde \alpha^{0\pi}_k - \alpha^{0\pi}_k$ satisfies
$$
J\cdot (\widetilde \alpha^0_k - \alpha_0^k) = - \sqrt{-1} (\widetilde \alpha^{0\pi}_k -\alpha^{0\pi}_k)
$$
and so is of type $(0,1)$. 

Finally we have
$$
\widetilde \Theta^0 = d\lambda + \sum_k \widetilde \alpha_k^0 \wedge \theta^k 
+  \widetilde \alpha^0 \wedge \lambda
$$
from the first of the structure equation \eqref{eq:Theta0i}, and $\widetilde \alpha_0^0 = 0$.
Therefore
$$
\widetilde \Theta^0 - \Theta^0 = \sum_k (\widetilde \alpha_k^0 -\alpha_k^0) \wedge \theta^k
+ (\widetilde \alpha^0 - \alpha^0) \wedge \lambda.
$$
This implies 
\be\label{eq:(1,1)-component}
(\widetilde \Theta^0 - \Theta^0)^{\pi(1,1)} = \sum_k (\widetilde \alpha_k^0 -\alpha_k^0)^{(0,1)} \wedge \theta^k.
\ee
Since $(\widetilde \Theta^0)^{\pi(1,1)} = 0 = (\Theta^0)^{\pi(1,1)} $ by Axiom (6), 
by taking the $(1,1)$-component of this equation,  we have derived
$$
\sum_k (\widetilde \alpha_k^{0} - \alpha_k^{0})^{(0,1)} \wedge \theta^k = 0.
$$
Since $\{\theta^k\}$ are linearly independent and of type $(1,0)$ and 
$(\widetilde \alpha_k^{0} - \alpha_k^{0})^{(0,1)}$ is of type
$(0,1)$, we derive  
\be\label{eq:alphak0(0,1)}
(\widetilde \alpha_k^{0})^{(0,1)}= (\alpha_k^{0})^{(0,1)}
\ee
for all $k$.  By taking the complex conjugate of \eqref{eq:(1,1)-component} and applying the above
argument thereto, we also derive
\be\label{eq:alphak0(1,0)}
(\widetilde \alpha_k^{0})^{(1,0)}= (\alpha_k^{0})^{(1,0)}.
\ee
For the vertical component, we have
$$
\widetilde \alpha_k^{0\perp} = \widetilde \alpha_k^0(R_\lambda)\, \lambda
$$
by definition. We evaluate the structure equation \eqref{eq:torsion-Theta0} against $(\eta_i,R_\lambda)$
applied to $\widetilde \Theta$, and obtain
$$
\widetilde \Theta^0(\eta_i,R_\lambda) = -\widetilde \alpha_i^0(R_\lambda)
$$
using $\widetilde \alpha_0^0 = 0$. Since $R_\lambda \intprod \widetilde \Theta = 0$, this implies
$$
\widetilde \alpha_i^0(R_\lambda) = 0
$$
which in turn implies 
$\widetilde \alpha_k^{0\perp} = 0 = \alpha_k^{0\perp}$. Combining  this with 
\eqref{eq:alphak0(0,1)} and \eqref{eq:alphak0(1,0)},  we have
proved
\be\label{eq:alpha-equal}
\widetilde \alpha_k^0 = \alpha_k^0 \quad \text{\rm for all $k$}.
\ee
On the other hand, we have the $(\pi 0)$-component 
$$
(\widetilde \Theta^i - \Theta^i)^{\pi 0} = \sum_i (\widetilde \beta^i - \beta^i) \wedge \lambda.
$$
This vanishes by the Axiom $R_\lambda\intprod T = 0$.  On the other hand, the skew-Hermitian
property, we have
$$
\beta^i(R_\lambda) = - \alpha^0_i(R_\lambda) = 0.
$$
Then Lemma \ref{lem:Axiom(2)(4)} implies
\be\label{eq:beta-equal}
\widetilde \beta^i = \beta^i.
\ee
Combining  \eqref{eq:omega-equal},  \eqref{eq:alpha-equal} and \eqref{eq:beta-equal}, 
we have proven $\widetilde \Theta = \Theta$ and
hence the uniqueness.

\subsection{Existence of contact triad connection}

We will divide our construction into the tasks of constructing $\omega = (\omega_j^i)$, and
$$
\alpha^0 = (\alpha_k^0), \quad \beta_0^t = (\beta_0^i) = (\beta^i),
$$
respectively.

Let $\widetilde \Omega$ be the connection one-form of
\emph{any} unitary connection of the complex vector bundle $\xi \otimes \C$ 
and its extension $\widetilde \Omega$
to $(\xi \otimes \C) \oplus \C \langle R_\lambda \rangle$. For example, we start with 
the matrix $\widetilde \Omega$ associated to the Levi-Civita connection.  Then
using the fact that any other affine connection $\nabla$ can be determined by
writing down the \emph{tensor}
$$
\nabla - \widetilde \nabla=: B
$$
Starting from the Levi-Civita connection and this fact, 
we will construct a new connection $\nabla$ by identifying its associated connection matrix 
$(\omega_i^j)$, $(\alpha_0^j), \, (\beta^0_j)$ of the given
almost contact frame.

We  express the Levi-Civita connection by the matrix
$$
\widetilde \Omega = \left(\begin{matrix} \widetilde \omega_i^j & \widetilde \alpha_k^{0\pi}\\
(\widetilde \beta^{k\pi})^t & 0 \end{matrix} \right)
$$
with respect to an almost contact frame as before.  We first prove the following analogue  of the Lemma \ref{lem:Axiom(2)(4)} for the Levi-Civita connection.

\begin{lem}\label{lem:Levi-Civita} Let $\widetilde \nabla$ be the Levi-Civita connection of the triad metric of
$(M, \lambda, J)$. Then
\bea
\widetilde \alpha_0^0 &=& 0 \,\, = \widetilde \beta^0, \label{eq:tildealpha00}\\
 \widetilde \beta^k & = & \overline{\widetilde{\alpha}^0_k} \quad \text{\rm for }\, k = 1, \cdots, n. \label{eq:tildealphabar}
 \eea
\end{lem}
\begin{proof} We have $\langle R_\lambda, Y \rangle = 0$ for any
vector field $Y \in \xi$. By differentiating the equation with respect to $R_\lambda$, we obtain
\beastar
\langle \widetilde \nabla_{R_\lambda} R_\lambda, Y \rangle & = & 
- \langle R_\lambda, \widetilde \nabla_{R_\lambda} Y \rangle\\
& = & - \langle R_\lambda, \widetilde \nabla_Y R_\lambda + [R_\lambda,Y] \rangle
\eeastar
where the second equality follows since the Levi-Civita connection is torsion free.
We compute
$$
\langle R_\lambda, \widetilde \nabla_Y R_\lambda + [R_\lambda,Y] \rangle
 = \frac12 Y \langle R_\lambda, R_\lambda \rangle + \langle R_\lambda, [R_\lambda,Y] \rangle = 0:
$$
Here the second equality follows since 
$\langle R_\lambda, R_\lambda \rangle \equiv 1$ and $\CL_{R_\lambda}$ preserves $\xi$ 
and so $[R_\lambda,Y] = \CL_{R_\lambda} Y \in \xi$.

Combining the above two equations, we have derived
$$
\langle \widetilde \nabla_{R_\lambda} R_\lambda, Y \rangle = 0
$$
for all $Y \in \xi$. On the other hand, we also have 
$$
\langle \widetilde \nabla_{R_\lambda} R_\lambda, R_\lambda \rangle 
= \frac12 R_\lambda \langle  R_\lambda, R_\lambda \rangle = 0.
$$
This proves $\widetilde \nabla_{R_\lambda} R_\lambda = 0$
which is equivalent to $\widetilde \alpha_0^0 = 0$.
Now both $\widetilde \beta^0 = 0$ and \eqref{eq:tildealphabar} follow from the 
skew-Hermitian property of the relevant connection matrix
is finishes the proof.
\end{proof}

Starting from this Levi-Civita connection, 
we will construct a new connection $(\omega_i^j)$, $(\alpha_0^j), \, (\beta^i)$
by considering the first structure equations \eqref{eq:torsion-Thetai},
\eqref{eq:torsion-Theta0} which are universal for general Riemannian connection.
Motivated by  Lemma \ref{lem:Levi-Civita}, we start with putting
\be\label{eq:alpha00-beta0k}
\alpha_0^0=0
\ee
which will be shown to be consistent with the first structure equations 
of the constructed connection at the end of the construction process.

Then we will determine the candidate connection so that associated torsion form satisfies
the following:
\begin{enumerate}
\item $\Theta$ satisfies \eqref{eq:Theta0i}.
\item  $\Theta^{\pi(1,1)} = 0$ and $\Theta^0 = d\lambda + \sum_k \alpha_k^{0\pi} \wedge \theta^k$.
\item $\alpha_k^{0\perp}= 0$ and $J\cdot \alpha_k^{0\pi} = -\sqrt{-1} \alpha_k^{0\pi} + c\, \theta^k$.
\end{enumerate}

We start with determining $\omega$. For this,  we will adapt
 the argument used in the proof of Theorem 2.1 \cite{kobayashi:canonical}
to the current contact triad context.
For this purpose, we first express
the $(1,1)$-component of $\widetilde \Theta^\pi$ of the connection $\widetilde \Omega$ into
\be\label{eq:tilde(1,1)}
(\widetilde \Theta^j)^{(1,1)} = \sum A_{i\bar k}^j   \theta^i \wedge \overline \theta^k 
=
 \sum (- A_{i\bar k}^j)  \overline \theta^k \wedge \theta^i
\ee
for some set of complex-valued functions $A_{i\bar k}^j$. Since $\theta^k \wedge \theta^i 
= - \theta^i \wedge\theta^k$, \emph{without loss of any generality, we may  assume}
\be\label{eq:A-skewsymmetry}
A_{i\bar k}^j = - A_{k \bar i}^j
\ee
after some rearrangement of the summands of \eqref{eq:tilde(1,1)}.
Then we define the new connection forms $\omega_i^j$ by putting
$$
\omega_i^j = \widetilde \omega_i^j + \sum_k  \overline A_{i \bar k}^j \theta^k +  A_{i\bar k}^j    \overline \theta^k.
$$
Then using the first structure equation \eqref{eq:Theta0i}, we compute
\bea\label{eq:Thetaj}
\Theta^j & = & d\theta^j + \sum_i \omega^j_i \wedge \theta^i \nonumber \\
& = & \left(\sum_i \widetilde \omega_i^j \wedge \theta^i + \widetilde \Theta^j \right) 
+ \left( \sum_i \widetilde \omega_i^j \wedge \theta^i + \sum_{i,k} \overline A_{i \bar k}^j  \theta^k \wedge \theta^i
+ \sum_{i,k}  A_{i \bar k}^j  \overline \theta^k \wedge \theta^i\right) \nonumber \\
& = &  \widetilde \Theta^j  + \sum_{i,k}  \overline A_{i \bar k}^j \theta^k \wedge \theta^i
+ \sum_{i,k} A_{i \bar k}^j  \overline \theta^k \wedge \theta.
\eea
From this, by taking the $(1,1)$-component of the equation, we obtain
\beastar
(\Theta^j)^{(1,1)} & =  & (\widetilde \Theta^j)^{(1,1)} + \sum_k A_{i \bar k}^j  \overline \theta^k \wedge \theta^i \\ 
& = &  \sum_k  - A_{i \bar k}^j   \overline \theta^k \wedge \theta^i 
+ \sum_k \ A_{i \bar k}^j   \overline \theta^k \wedge \ \theta^i
 \\
& = &  \sum_k (- A_{i \bar k}^j +  A_{i\bar k}^j) \overline \theta^k \wedge \theta^i = 0
\eeastar
where the second equality follows from \eqref{eq:tilde(1,1)}.

Then we can write the $0$-component  $\widetilde \Theta^0$ of the first structure equation of the Levi-Civita connection
$\widetilde \Omega$ into
$$
\widetilde \Theta^0 = d\lambda + \sum_k \widetilde \alpha_k^{0\pi} \wedge \theta^k.
$$
By decomposing $\widetilde \alpha_k^{0\pi}$ into 
$\widetilde \alpha_k^{0\pi} = (\widetilde \alpha_k^{0\pi})^{(0,1)} + (\widetilde \alpha_k^{0\pi})^{(1,1)}$,
we obtain
\beastar
(\widetilde \Theta^0)^{(1,1)} & = & d\lambda + \sum_k (\widetilde \alpha_k^{0\pi})^{(0,1)} \wedge \theta^k\\
(\widetilde \Theta^0)^{(2,0)} & = & \sum_k (\widetilde \alpha_k^{0\pi})^{(1,0)} \wedge \theta^k.
\eeastar
Furthermore we obtain
$$
J \cdot \widetilde \alpha_k^{0\pi} =  \sqrt{-1}(\widetilde \alpha_k^{0\pi})^{(0,1)} - \sqrt{-1} (\widetilde \alpha_k^{0\pi})^{(0,1)}
$$
while
$$
-\sqrt{-1} \cdot \widetilde \alpha_k^{0\pi} = -\sqrt{-1} (\widetilde \alpha_k^{0\pi})^{(0,1)} - \sqrt{-1} (\widetilde \alpha_k^{0\pi})^{(0,1)}.
$$
Hence
$$
J \cdot \widetilde \alpha_k^{0\pi} = - \sqrt{-1} \cdot \widetilde \alpha_k^{0\pi} + 2\sqrt{-1} (\widetilde \alpha_k^{0\pi})^{(1,0)}.
$$

Now to determine $\alpha_k^{0\pi}$, we first put
$$
(\alpha^0_k)^\perp : =( \widetilde \alpha^0_k)^\perp = 0
$$
using the vanishing of 
$\alpha_k^{0\perp} = \alpha_k(R_\lambda) \lambda = 0$. Then 
we decompose 
$$
\alpha_k^0 = \alpha_k^{0\pi} = (\alpha_k^{0\pi})^{(1,0)} + (\alpha_k^{0\pi})^{(0,1)}
$$
and put
\be\label{eq:new-alpha}
\begin{cases}
 (\alpha_k^{0\pi})^{(0,1)}: =  (\widetilde \alpha_k^{0\pi})^{(0,1)}, \\
 (\alpha_k^{0\pi})^{(1,0)}: = J \cdot  \left(\sum_\ell B_{k\ell} \theta^\ell \right) + 
\sqrt{-1}\left ( \sum_\ell B_{k\ell} \theta^\ell\right)
 \end{cases}
 \ee
for some coefficients $B_{k\ell}$, and solve for $B_{k\ell}$ so that $\alpha_k^0$ satisfies Condition (3) above.  
We compute
\beastar
&{}& 
J \cdot \alpha_k^{0\pi} + \sqrt{-1} \cdot \alpha_k^{0\pi} \\
& = & \left(J (\alpha_k^{0\pi})^{(0,1)} + \sqrt{-1} (\alpha_k^{0\pi})^{(0,1)}\right) +
\left(J \cdot (\alpha_k^{0\pi})^{(1,0)}+ \sqrt{-1} ( \alpha_k^{0\pi})^{(1,0)}\right) \\
& = &  \left(J (\widetilde \alpha_k^{0\pi})^{(0,1)} + \sqrt{-1} (\widetilde \alpha_k^{0\pi})^{(0,1)}\right) 
+ J \cdot  \left(\sum_\ell B_{k\ell} \theta^\ell \right) + 
\sqrt{-1}\left ( \sum_\ell B_{k\ell} \theta^\ell\right) \\
& = &  \left(- \sqrt{-1} (\widetilde \alpha_k^{0\pi})^{(0,1)} + \sqrt{-1} (\widetilde \alpha_k^{0\pi})^{(0,1)}\right) 
+   \left(\sum_\ell B_{k\ell} \sqrt{-1} \theta^\ell \right) + 
\sqrt{-1}\left ( \sum_\ell B_{k\ell} \theta^\ell\right) \\
& = &  2 \sqrt{-1}\left ( \sum_\ell B_{k\ell} \theta^\ell\right).
\eeastar
Therefore we need $B_{k\ell} $ to satisfy the equation
$$
2 \sqrt{-1}\left ( \sum_\ell B_{k\ell} \theta^\ell\right) = c\,  \theta^k.
$$
By solving this equation, we get the unique solution
$$
B_{k\ell} = \frac{c}{2\sqrt{-1}} \delta_{k\ell}
$$
and hence we obtain
$$
\alpha_k^{0\pi} =  (\widetilde \alpha_k^{0\pi})^{(0,1)} + \frac{c}{2\sqrt{-1}} \theta^k.
$$
We  next determine $\beta^k_0 = \beta^k$ by the formula in Lemma \ref{lem:Axiom(2)(4)}
by putting $\beta^{k\pi}: = -\overline \alpha_k^{0\pi}$ and
 $(\beta^k)^\perp = 0$  which are forced by
 the skew-Hermtian property of the connection matrix $\Omega$ and by $(\alpha^0_k)^\perp = 0$.
 
Finally it remains to check consistency of the  choice $\Omega$ of the matrix-valued 1-form 
made as above by checking the following:
\begin{enumerate}
\item The above constructed
 $\Omega$ defines a connection form with respect to the given (dual) almost contact frame
$\{\theta_1, \cdots, \theta_k\} \cup \{\lambda\}$ on $\xi^{(1,0)} \oplus \C \langle \lambda \rangle$.
\item Then the form $\Theta$ satisfies the first structure equation of $\Omega$ so that
it is indeed the torsion form of $\Omega$.
\item 
Furthermore it also satisfies the defining property, Condition (3) above by construction, 
for a  $c$-contact triad connection  with respect to the given almost contact frame.
\end{enumerate}
However these requirements automatically hold by construction because they are precisely
the conditions that we have used to solve the above  matrix elements $(\omega^i_j),\, (\alpha_k^0)$ and
$(\beta^k_0)$.
This establishes the existence of a required triad connection and finish the proof.
\qed

\subsection{Derivation of additional properties of the connection}

Another consequence of the above proof, \eqref{eq:dlambda}
 and the expression \eqref{eq:Thetaj} of the connection form $\Theta = (\Theta^j)$
is the following.

\begin{thm}\label{thm:(2,0)-vanishing}
 Let $(M,\lambda,J)$ be any contact triad. Then for any triad connection $\nabla$, its connection forms 
$\Theta^j$ for $j = 1, \cdots, n$ is of $(0,2)$-type.
\end{thm}
\begin{proof} By taking the differential of \eqref{eq:dlambda}, we obtain
$$
0 = \sum_k d\theta^k \wedge \overline \theta^k + \sum_k \theta^k \wedge d\overline \theta^k.
$$
Substituting the first structure equation $d\theta^k = \Theta^k - \sum_{\ell =1}^n \omega^k_\ell \wedge \theta^\ell$ hereinto, we 
obtain
\bea\label{eq:0}
0 & = & \sum_k \sum_\ell \left(\Theta^k - \sum_{\ell =1}^n \omega^k_\ell \wedge \theta^\ell\right) \wedge \overline \theta^k 
+ \sum_k \theta^k \wedge \left(\overline \Theta^k - \sum_{\ell =1}^n \overline \omega^k_\ell \wedge \overline \theta^\ell\right) 
\nonumber\\
& = & \sum_k \sum_\ell \Theta^k  \wedge \overline \theta^k 
+ \sum_k \theta^k \wedge \overline \Theta^k.
\eea
Using the vanishing of $(\Theta^k)^{(1,1)} = 0$, the $(2,1)$-component of the right hand side becomes
$$
\sum_k \theta^k \wedge (\overline \Theta^k)^{(2,0)}.
$$
On the other hand, we also derive
$$
(\Theta^k)^{(2,0)} = \sum_{i,\ell}  \overline A_{i \bar \ell }^k \theta^\ell \wedge \theta^i
= \sum_{i > \ell} ( \overline A_{i \bar \ell }^k  - \overline A_{\ell \bar i}^k)\theta^\ell \wedge \theta^i
= 2 \sum_{i > \ell} \overline A_{i \bar \ell }^k \, \theta^\ell \wedge \theta^i
$$
from \eqref{eq:tilde(1,1)}, \eqref{eq:A-skewsymmetry} and \eqref{eq:Thetaj}. Therefore by taking $(2,1)$-component of \eqref{eq:0}, we have obtained
\beastar
0 & = & \sum_k \theta^k \wedge \overline{\left(2 \sum_{i > \ell} \overline A_{i \bar \ell }^k\, \theta^\ell \wedge \theta^i\right)}
= 2 \sum_k \sum_{i>\ell} A_{i \bar \ell}^k \theta^k \wedge \overline \theta^\ell \wedge \overline \theta^i
\eeastar
which is equivalent to $A_{i\bar \ell}^k= 0$ for all $i, \, k, \,\ell$. This proves the corollary.
\end{proof}

\begin{rem}\label{rem:ac-extension}
\begin{enumerate}
\item 
The proof of Theorem \ref{thm:(2,0)-vanishing} follows the same reasoning of 
that of Theorem \ref{thm:aK-(2,0)-vanishing} from \cite{kobayashi:canonical}: 
The proof is based on the \emph{closedness of its fundamental 2-form
$\omega: = g(J\cdot, \cdot)$} of any almost Hermitian metric. (See \cite[Theorem 7.1.14]{oh:book1} too.)
This corresponds to the closedness of $d\lambda = g(J\cdot, \cdot)$
for the contact triad metric $g$ in the current context.
In this regard, our contact triad connection is the contact counterpart of 
Kobayashi's notion of canonical connection on almost K\"ahler manifolds \cite{kobayashi:canonical}.
\item In fact, our construction of ($c$-)contact triad connection can be verbatim extended to
the context of \emph{almost contact manifolds} \cite{gray} with the replacement of the triad $(M,\lambda,J)$ 
as follows.
Consider the almost contact manifold $(M,J, g)$ with the splitting
$$
TM = \xi \oplus \underline{\R}, \quad J^2 = -\Pi, \quad \xi \perp_g \underline{\R}
$$
where $\Pi: TM \to TM$ is the idempotent with $\Im \Pi = \xi, \quad \ker \Pi = \underline{\R}$.
We define the two-form
\be\label{eq:Phi-contact}
\Phi: = g(J\cdot, \cdot)
\ee
which we call the \emph{fundamental two form} of the almost contact manifold $(M,J,g)$
 following the terminology used for the two-form \eqref{eq:Phi} of almost Hermitian
manifolds in \cite{kobayashi-nomizu2}.
Then all the discussion go through except Theorem \ref{thm:(2,0)-vanishing}, \emph{unless $\Phi$ is a
closed two-form}. (Such an almost contact structure is called a \emph{quasi-contact structure} in 
\cite{casals-pancho-presas}.) When $\Phi$ is exact,  say $\Phi = d\lambda_{(J,g)}$, then $\lambda_{(J,g)}$
defines a contact form.
\end{enumerate}
\end{rem}

\end{document}